\documentclass[preprint,12pt]{elsarticle}

\usepackage[tbtags]{amsmath}
\usepackage{bm}
\usepackage{amsfonts}
\usepackage{graphics}
\usepackage{graphicx}
\usepackage{psfrag}
\usepackage{eurosym}
\usepackage{pstricks}
\usepackage{array}
\usepackage{shortvrb}
\usepackage{epsf}
\usepackage{graphicx}
\usepackage{rotating}
\usepackage{float}
\usepackage{color}
\usepackage{multirow, makecell}
\usepackage{float}
\usepackage{colortbl}
\usepackage{ifpdf}
\usepackage{multicol}

\usepackage{hyperref}
\hypersetup{colorlinks=true, linkcolor=blue, anchorcolor=red, citecolor=blue, filecolor=red, urlcolor=red, pdfauthor=author}
            
\usepackage{subfigure}

\usepackage{geometry}
\usepackage{alphalph}
\usepackage{etoolbox}
\usepackage{pdfpages}

\usepackage[T1]{fontenc}

\usepackage{caption}

\newcommand{\comment}[1]{}

\patchcmd{\subequations}{\alph{equation}}{\alphalph{\value{equation}}}{}{}
\def \Noin {{\hskip 3pt \rm /\kern -9pt \in\hskip 1pt}}

\def \R {{\rm I\kern -2.2pt R\hskip 1pt}}

\allowdisplaybreaks

\begin{document}

\begin{frontmatter}

\title{A Neural Network-Based Distributional Constraint Learning Methodology for Mixed-Integer Stochastic Optimization}

\author[inst1]{Antonio Alc\'antara\corref{cor1}}
\ead{antalcan@est-econ.uc3m.es}
\author[inst1,inst2]{Carlos Ruiz}
\ead{caruizm@est-econ.uc3m.es}
\cortext[cor1]{Corresponding author}

\affiliation[inst1]{organization={Department of Statistics, University Carlos III of Madrid},
            addressline={Avda. de la Universidad, 30}, 
            city={Legan\'es},
            postcode={28911}, 
            state={Madrid},
            country={Spain}}

\affiliation[inst2]{organization={UC3M-BS Institute for Financial Big Data (IFiBiD), University Carlos III of Madrid},
            addressline={Avda. de la Universidad, 30}, 
            city={Legan\'es},
            postcode={28911}, 
            state={Madrid},
            country={Spain}}

\begin{abstract}

The use of machine learning methods helps to improve decision making in different fields. In particular, the idea of bridging predictions (machine learning models) and prescriptions (optimization problems) is gaining attention within the scientific community. One of the main ideas to address this trade-off is the so-called Constraint Learning (CL) methodology, where the structures of the machine learning model can be treated as a set of constraints to be embedded within the optimization problem, establishing the relationship between a direct decision variable $x$ and a response variable $y$. However, most CL approaches have focused on making point predictions for a certain variable, not taking into account the statistical and external uncertainty faced in the modeling process. In this paper, we extend the CL methodology to deal with uncertainty in the response variable $y$. The novel Distributional Constraint Learning (DCL) methodology makes use of a piece-wise linearizable neural network-based model to estimate the parameters of the conditional distribution of $y$ (dependent on decisions $x$ and contextual information), which can be embedded within mixed-integer optimization problems. In particular, we formulate a stochastic optimization problem by sampling random values from the estimated distribution by using a linear set of constraints. In this sense, DCL combines both the high predictive performance of the neural network method and the possibility of generating scenarios to account for uncertainty within a tractable optimization model. The behavior of the proposed methodology is tested in a real-world problem in the context of electricity systems, where a Virtual Power Plant seeks to optimize its operation, subject to different forms of uncertainty, and with price-responsive consumers.
\end{abstract}

\begin{keyword}
Stochastic optimization \sep Constraint learning \sep Distribution Estimation \sep Neural Networks \sep Mixed-Integer Optimization
\end{keyword}

\end{frontmatter}

\section{Introduction}
\label{sec:intro}

Uncertainty exists in most areas of our daily lives, from weather conditions to energy supplies, and mathematical modeling is also an area that is not exempt from this uncertainty. In this sense, both predictive and prescriptive processes are exposed to different degrees of uncertainty that can greatly affect their performance.

In the case of predictive problems, probabilistic forecasting is one of the main approaches in the literature to deal with uncertainty, which estimates a probability for every possible event that may occur. This forecasting methodology can be used in a wide range of applications, such as estimating electricity price \cite{nowotarski2018recent}, photovoltaic generation \cite{wan2016probabilistic}, or even COVID-19 mortality \cite{taylor2021combining}. On the other hand, standard point prediction techniques assign a single value to a target variable (typically the average, as in linear least squares, or a label, in classification), making this methodology not adequate for addressing uncertainty in the predictions. 

The state-of-the-art in probabilistic forecasting is composed of regression-based methods that estimate Prediction Intervals (PIs), i.e., lower and upper bounds that contain a specific target variable with a certain probability. One of the main approaches is to estimate quantiles that act as PIs, fitted by minimizing quantile loss \cite{hao2007quantile}, or setting a specific coverage to obtain bounds for PIs, which is done primarily by neural networks (NNs). This field gained relevance with the introduction of the Coverage-Width Criterion (CWC), used as a training loss and based on the quality of estimated PIs \cite{khosravi2010lower}. More recently, Quality-Driven (QD) loss was introduced \cite{pearce2018high}, solving some previous issues that can occur while training with CWC and improving PI quality.

However, when employing these aforementioned methods, the practitioner typically chooses a nominal coverage or probability, which disregards the tails of the prediction distribution,  which may have a strong impact on the decision-making process. One solution to incorporate this information is to employ predictive methods that allow us to estimate the complete probability distribution of the forecast. For example, Deep Autoregressive Recurrent Networks (DeepAR) \cite{salinas2020deepar} are designed to estimate the mean and variance of the possible outcome in a time series context.

When it comes to making decisions, i.e., prescriptions, decision-makers typically suffer from high levels of uncertainty within their optimization problems. The classical literature on optimization has been based on assuming a given probability distribution or uncertainty set for the different parameters that participate in the problem \cite{murty1994operations}. Nowadays, in many relevant applications, the wide availability of data allows us to make estimations of these parameters without needing to assume any distribution over them. In this regard, traditional stochastic optimization \cite{birge2011introduction} can be adapted to this empirical data-driven approach. In particular, recent research has been developed to study how to optimally work with predictions and prescriptions at the same time. From the classical ``Predict, then Optimize'' where parameters were estimated making use of statistical or machine learning models to be used later in optimization problems, new approaches have been derived to address the trade-off between predictive accuracy and prescriptive performance \cite{bertsimas2020predictive, elmachtoub2022smart}.

These aforementioned methodologies try to estimate parameters $y$ assuming independence with respect to the optimization decision variables $x$. However, there are relevant situations where $y$ is (in some sense) dependent on direct decisions $x$, which can be addressed by the field of Constraint Learning (CL) \cite{fajemisin2021optimization}. This methodology is characterized by explicitly modeling the relationship between $y$ and $x$ thanks to the use of different machine learning models, which can be embedded within mixed-integer optimization problems as a set of piece-wise linear constraints.

The CL paradigm has been gaining attention among researchers and has been employed in different fields. For example, a neural network model was embedded to learn power flow constraints in \cite{kody2022modeling}, while in \cite{grimstad2019relu} it was employed to solve a production optimization problem from an offshore petroleum production field. In \cite{maragno2021mixed}, a chemotherapy regimen design was planned using constraint learning for different toxicities and survival rates, while in \cite{mistry2021mixed} gradient-boosted trees were embedded into Mixed-Integer convex non-linear optimization problems. Finally, in \cite{bergman2022janos}, an end-to-end software is developed to train machine learning models and embed them into optimization problems.

Nevertheless, the point prediction machine learning models that are used in the previous CL examples may not be adequate in some contexts where statistical uncertainty is present. We can take as an example some demand $y$ that is elastic with respect to a market price $x$. In this case, the possible stochastic values (due to some statistical uncertainty) that demand can take given a particular price $x$ will directly influence the profit that a seller can make.

To address these issues, in this paper, we propose a methodology based on CL that allows us to incorporate uncertainty in the form of scenario realizations sampled from estimated conditional distributions: Distributional Constraint Learning (DCL). We make use of modern neural network structures for distribution forecasting that can be embedded as a set of constraints within mixed-integer stochastic optimization problems in order to address the statistical uncertainty of response decision variables.

The proposed methodology will be tested in the optimal operation of a Virtual Power Plant (VPP). A VPP works as an intermediary agent between the Market Operator (MO) and a collection of consumers and Distributed Energy Resources (DERs). The aim of the VPP is to gather small consumers and DERs and trade energy with the electricity market. In other words, the VPP manages the energy resources available in the DERs to meet the demand requirements of the participants and by setting the appropriate price-tariffs for the energy transactions.

Several studies have been published on VPP operations. For example, in \cite{babaei2019data}, wind power generation is included as DER, and a robust two-stage optimization problem is presented to solve the unit commitment. In \cite{zhao2021operating}, PIs are obtained to address the problem of calculating operating reserve. However, the structure of the VPP makes its day-ahead operating optimization problem suitable to benefit from the proposed stochastic DCL methodology. For instance, a realistic assumption motivated by current demand response programs would be that the elasticity demand of the participants in a VPP is price-sensitive, i.e., the price-tariff that the operator sets to the DERs will influence their demand. In this sense, our approach will allow the decision-maker to generate high-quality forecasting demand response scenarios, which are conditioned by its own price decisions. These functional relationships are embedded within its optimization problem, while accounting for other sources of uncertainty.

Therefore, we summarize the main contributions of this work concerning the state-of-the-art on the aforementioned topics as the following:

\begin{itemize}
    \item[-] to extent the constraint learning methodology in stochastic mixed-integer optimization by addressing the statistical uncertainty in the response variables.
    \item[-] to deal with this uncertainty by estimation conditional distributions from modern neural network-based methods, whose structure allows to be embedded within the optimization problem as sets of constraints, while maintaining high predictive performance and the possibility of generating stochastic scenarios. We denote this novel methodology as Distributional Constraint Learning (DCL).
    \item[-] to develop an open-access software tool ``DistCL'' that allows practitioners to employ the DCL methodology in their stochastic optimization problems.
    \item[-] to test and illustrate the validity of the approach in a real-world-based Virtual Power Plant optimization problem. 
\end{itemize}

The structure of this article is as follows. Section \ref{sec:cl_background} will address statistical uncertainty and constraint learning methodology in more detail. Section \ref{sec:distcl} will introduce the proposed DCL methodology, from the neural network training process to the embedding procedure and the scenario generation in the stochastic mixed-integer optimization problem. Then, the Virtual Power Plant case study will be presented and discussed in Section \ref{sec:case_studies}. Finally, Section \ref{sec:conclu} draws the main conclusions of this work.

\section{Background: Uncertainty in Constraint Learning}
\label{sec:cl_background}

The constraint learning methodology seeks to explicitly model the relationship between a direct decision variable $x$ and a response decision variable $y$. This means that the value of $y$ is directly dependent on that assigned to $x$ and, in general, can also be affected by exogenous contextual information $\theta$.

We can take as an example the optimization model (\ref{eq:cl_opt}), which aims to minimize the expected value of the cost function $c(x,y;\theta,\xi): \mathbb{R}^{d_x} \rightarrow \mathbb{R}$ (\ref{eq:cl_OB}). This cost depends not only on the decision variables $x \subset \mathbb{R}^{d_x}$ and $y \subset \mathbb{R}^{d_y}$, but also on external contextual information $\theta \subset \mathbb{R}^{d_{\theta}}$ and uncertainty $\xi \subset \mathbb{R}^{d_{\xi}}$. Note that additional constraints could be applied through known functions $g(\cdot)$ on $x$ and $y$, as stated in the constraint (\ref{eq:cl_cons1}).

\begin{subequations}\label{eq:cl_opt}
\begin{align}
\underset{x, y}{\min} \;& \mathbb{E}\left[ c(x,y;\theta,\xi) \right] \label{eq:cl_OB}\\
\text{s.t.}&\notag \\
  &g(x,y)\leq 0  \label{eq:cl_cons1}\\
  &y = f^{\mathcal{D}}(x;\theta,\xi) \label{eq:cl_cons2} 
\end{align}
\end{subequations}

The learned constraint is established in (\ref{eq:cl_cons2}). As can be seen, the relationship between $y$ and $x$ (and contextual information $\theta$ and uncertainty $\xi$) is modeled using a selected point prediction model $f^{\mathcal{D}}$, which is trained using a data set $\mathcal{D} = \{x_i,y_i,\theta_i \}_{i=1}^N$. 

Regarding this prediction model, if we consider (\ref{eq:cl_opt}) as a mixed-integer optimization problem, it can be selected from a wide variety of machine learning methods, such as linear models, or piece-wise linearizable decision tree-based methods or neural networks. All these methods have in common a structure that allows the decision-maker to embed them within the mixed-integer problem as a set of piece-wise linear constraints \cite{fajemisin2021optimization}. 

However, some concerns arise from the constraint learning methodology. First, the decision-maker may be biased due to poor modeling of $y$, i.e., a high fitting error from $f^{\mathcal{D}}(\cdot)$. In this sense, the training process plays an important role, as it is aimed at a minimum predictive error. There are even cases where the decision-maker may need probabilistic guarantees over the learned constraint, and the employment of chance constraints is required in the optimization problem. In these cases, as classical point estimation methods are not appropriated, quantile estimation-based methods are used to characterize these constraints \cite{alcantara2022data}.

Another topic of concern is the uncertainty associated with the learned variable $y$. From a machine learning perspective, there is a systematic uncertainty added to the problem when the model $f^{\mathcal{D}}(\cdot)$ is employed, which we refer to as \textit{uncertainty in the weights}. For example, when using complex deep neural networks as a predictive model, different initializations in the training process may lead to slightly different results due to small changes in model weights. In the case of a classification task, this could mean changing a prediction from one class to another in observations where the decision is not clear for the machine learning model. In the case of a regression task, additional constraints such as an upper or lower bound for $y$ in (\ref{eq:cl_opt}), could be fulfilled or not depending on the different weight realizations of the predictive model. 

On the other hand, uncertainty may be naturally present in the response variable $y$. This is what is called statistical uncertainty, and can be defined as the variability of the output $y$ given similar values of $x$. In this case, point estimations may not be appropriated to capture the uncertainty, and probabilistic forecasts are in order, for example, with PIs.

\begin{figure}[!ht]
    \centering
    \subfigure[]{\includegraphics[width=0.49\textwidth]{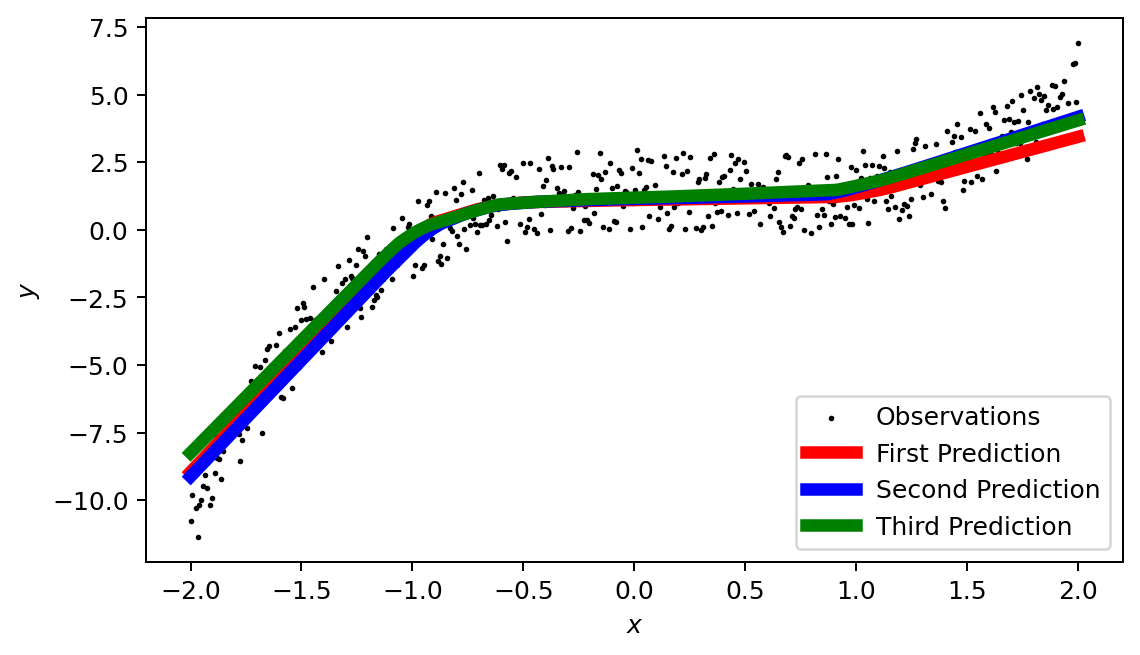}}
    \subfigure[]{\includegraphics[width=0.49\textwidth]{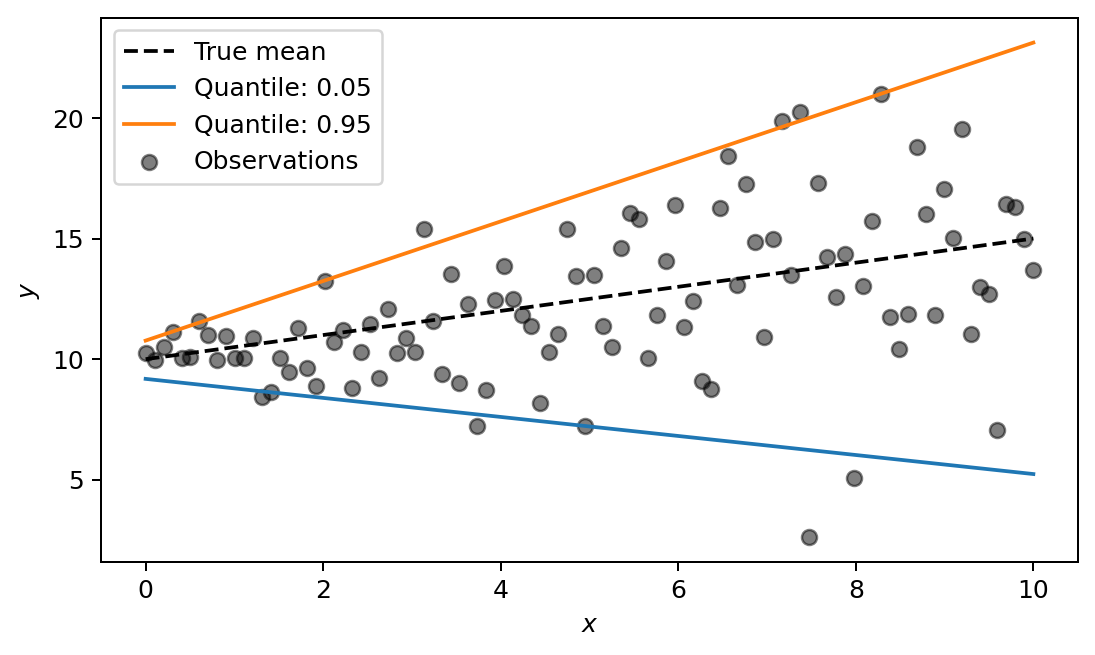}}
    \caption{Systematic (a) and statistical (b) uncertainty in machine learning}
    \label{fig:uncertainties}
\end{figure}

Graphically, the difference between systematic and statistical uncertainty can be seen in Figure \ref{fig:uncertainties}. On the left side (Figure \ref{fig:uncertainties} (a)), systematic weight uncertainty is represented. Note that different initializations may lead to slightly different predictions. On the right side (Figure \ref{fig:uncertainties} (b)), the statistical uncertainty is addressed using a linear quantile estimation method. We can notice how both quantiles represent the lower and upper values of a 90\% probability PI.

Returning to optimization problems with learned variables or constraints, the decision-maker can take into account systematic uncertainty through the use of ensemble methods, i.e., different predictions from different initialized models are ``ensembled'' to throw a combined prediction reducing systematic uncertainty.

However, statistical uncertainty can be addressed with probabilistic techniques. An important approach is to estimate a PI for the learned variable $y$ with a piece-wise linearizable machine learning method. Therefore, the uncertain value of $y$ would be contained between these two bounds. However, this PI-based approach would imply two main limitations:

\begin{enumerate}
    \item The decision-maker will have to decide on a coverage for the PI in advance, which will lead to a trimmed distribution. This may create a bias in the optimization problem when not taking into account the worst- or best-case scenarios.
    \item In the context of stochastic optimization, the generation of scenarios may not be trivial. For example, if we assume that the stochastic values of the learned variable $y$ must be contained in the PI, no additional information about the distribution of these stochastic values is available to sample them.
\end{enumerate}

Therefore, a relevant approach to take into account the issues mentioned above is to learn the complete conditional distribution of the learned variables. In this work, we introduce a new Distributional Constraint Learning (DCL) methodology that models the relationship between the distribution of $y$ and the decision variables $x$ with piece-wise linearizable machine learning models. These can be easily embedded in mixed-integer optimization problems as sets of constraints and adapted for stochastic optimization.

\section{Distributional Constraint Learning for Stochastic Optimization}
\label{sec:distcl}

As stated in Section \ref{sec:cl_background}, the Distributional Constraint Learning (DCL) methodology will be based on the estimation of the distribution of the response variable. In this section, we will introduce a neural network-based approach as one method to deal with the statistical uncertainty of $y$ through the output of its conditional distribution, which can be directly embedded within a mixed-integer optimization problem as a set of constraints. This methodology will allow practitioners to estimate the possible values that the response variable $y$ (dependent on decisions $x$) can take when decision-making takes place. The DCL framework can be employed as a useful tool to address uncertainty in the constraint learning methodology, which can have a strong impact on the prescription process.

\subsection{Optimization framework}
\label{sec:opt_dcl}

We define DCL as an extension of the constraint learning methodology to address uncertainty. An example of the DCL optimization framework is shown in (\ref{eq:dcl}) for a stochastic optimization problem.

\begin{subequations}\label{eq:dcl}
\begin{align}
\underset{x, \tilde{y}}{\min} \;& \mathbb{E}_{\tilde{y} \sim \mathcal{Y}(\hat{\mu}_y, \hat{\sigma}_y)}\left[ c(x,\tilde{y};\theta,\xi) \right] \label{eq:dcl_obj}\\
\text{s.t.}&\notag \\
  &g(x,\tilde{y})\leq 0  \label{eq:dcl_othercons}\\
  &\hat{\mu}_y, \hat{\sigma}_y = f^{\mathcal{D}}(x;\theta,\xi) \label{eq:dcl_modelcons}
\end{align}
\end{subequations}

We can see in (\ref{eq:dcl_modelcons}) how the estimations for the mean and standard deviation of the distribution of $y$ ($\hat{\mu}_y$ and $\hat{\sigma}_y$) are obtained from a piece-wise linearizable machine learning model $f^{\mathcal{D}}(\cdot)$, whose output is conditioned by decisions made over $x$. In this stochastic context, the random values $\tilde{y}$ come from the real unknown distribution of $y$ ($\mathcal{Y}$) with an estimated mean $\hat{\mu}_y$ and standard deviation $\hat{\sigma}_y$, as emphasized in the expected value calculated in the objective function (\ref{eq:dcl_obj}). However, note that this expected value may also account for other exogenous sources of uncertainty ($\xi$).

In this setup, it is important to correctly model the distribution of the learned variable $y$ within the optimization problem, as it may directly influence the cost function (\ref{eq:dcl_obj}): for example, when a price or a demand is learned from direct decisions and contextual information.

In the context of stochastic programming, we propose to generate scenario realizations for $\tilde{y}$ (i.e., $y_{\omega}$) directly from constraints within the optimization problem (\ref{eq:dcl}) by making some soft assumptions about the distribution $\mathcal{Y}$. For example, based on the Central Limit Theorem we can assume that the response distribution $\mathcal{Y}(\hat{\mu}_y, \hat{\sigma}_y)$ is actually a normal distribution $\mathcal{N}(\hat{\mu}_y, \hat{\sigma}_y)$. Therefore, the values for $\tilde{y}$ can be easily sampled using values from a standard normal distribution, as will be shown in Section \ref{sec:sce_gen}.

Furthermore, depending on the application, it may be suitable to consider $x$ as a first-stage decision variable, while each realization of the learned variable $y_{\omega}$ can be treated as a second-stage decision variable, but conditioned on the first-stage variable. Indeed, (\ref{eq:dcl}) can also be adapted to include other first- and second-stage decision variables not directly dependent on $x$, as well as other exogenous uncertainty sources that characterize each scenario.

\subsection{Distributional Learning with Neural Networks}
\label{sec:dnn}

One of the steps for the DCL methodology is to embed a piece-wise linearizable machine learning model within the optimization problem, which must be able to estimate a conditional mean and a standard deviation for the learned variable $y$. For this purpose, we present an already well-established neural network-based method for distributional estimation, which is denoted as Distributional Neural Network (DNN) \cite{nix1994estimating}.

\begin{figure}[!ht]
    \centering
    \includegraphics[width=\textwidth]{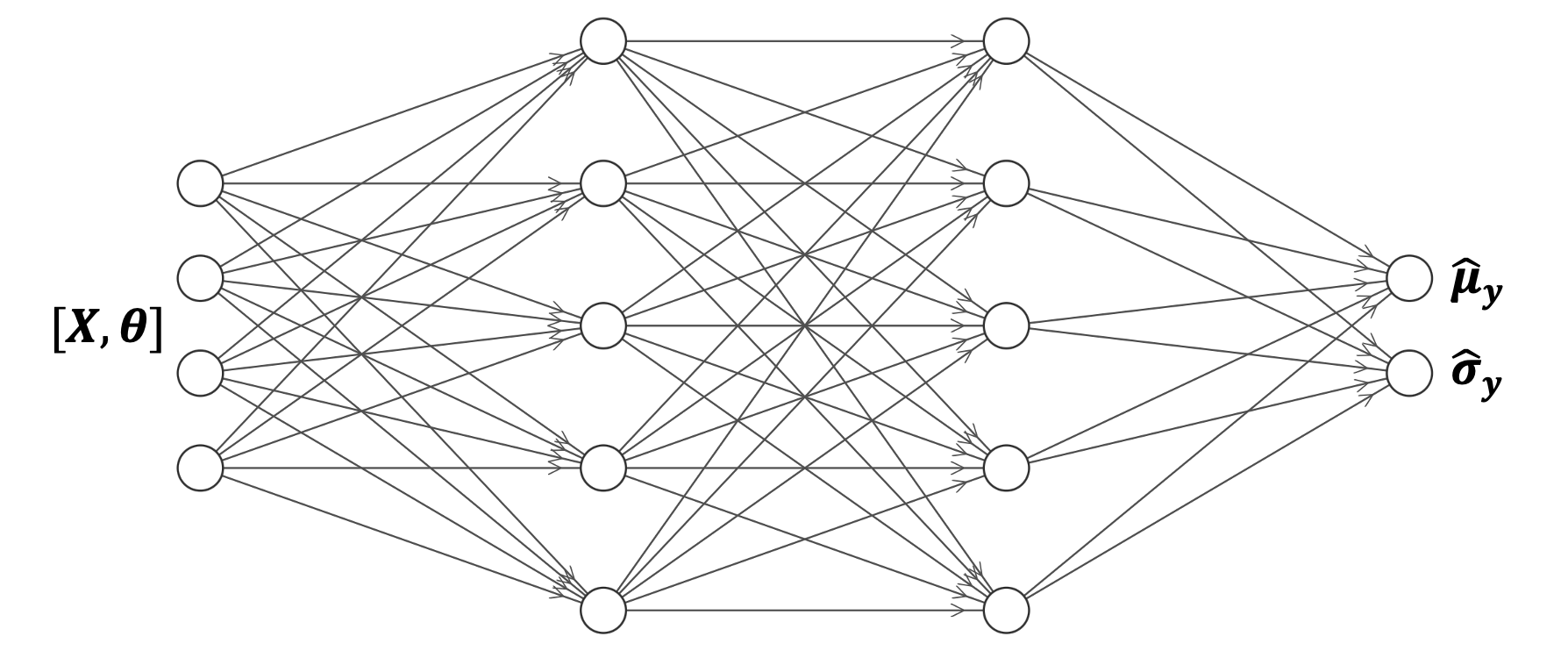}
    \caption{Distributional Neural Network structure}
    \label{fig:dnn_str}
\end{figure}

The structure of the DNN is the same as that of a classical feed-forward NN. Typically, they are made up of an input layer, $L-2$ hidden layers, and an output layer, with several neurons in each layer. A visual example is presented in Figure \ref{fig:dnn_str} for a DNN with 2 hidden layers. We can notice how the DNN takes as inputs previous decisions for $x$ and contextual information $\theta$ presented in the previously introduced dataset $\mathcal{D}$. As output, the DNN produces an estimate for the conditional mean and standard deviation of $y$, i.e., $\hat{\mu}_y$ and $\hat{\sigma}_y$, respectively.

In general, for a given hidden layer $l \in L$ of the DNN with nodes $N^l$, the value of a node $i \in N^l$, denoted as $v_i^l$, is calculated using the weighted sum of the node values of the previous layer, followed by a non-linear activation function $h(\cdot)$. This value is given as:
\begin{equation}
\label{eq:output_1ayer}  
 v_i^{l}=h\left(b_i^l + \sum_{j\in N^{l-1}} w_{ij}^{l} v_j^{l-1} \right)
\end{equation}

\noindent where $b_i^l$ represents the bias term, and $w_{ij}^{l}$ the weights for node $i$ in layer $l$. This computation will continue in the following layers, taking as inputs the outputs of the previous layers, until the final one is reached.

The DNN training process is based on the Gaussian negative log likelihood loss \cite{nix1994estimating}, assuming that the real distribution of $y$ ($\mathcal{Y}$) is normal or convergent to normal. The loss is defined in (\ref{eq:gaussianloss}).

\begin{equation}
\label{eq:gaussianloss}
\mathcal{L} = \frac{1}{2}\left( \log{\hat{\sigma}_y^2} + \frac{(y - \hat{\mu}_y)^2 }{\hat{\sigma}_y^2} \right)
\end{equation}

In practice, we fit the DNN with dataset $\mathcal{D} = \{x_i,y_i,\theta_i \}_{i=1}^N$. Therefore, the empirical Gaussian negative log likelihood loss is defined as a mean over losses in the dataset:

\begin{equation}
\label{eq:gaussianloss2}
\mathcal{L}_{\mathcal{D}} = \frac{1}{2N} \sum_{i=1}^N \left( \log{\hat{\sigma}_{y_i}^2} + \frac{(y_i - \hat{\mu}_{y_i})^2 }{\hat{\sigma}_{y_i}^2} \right)
\end{equation}

Note that $\hat{\sigma}_{y_i}$ and $\hat{\mu}_{y_i}$ are outputs of the DNN, which takes as inputs $x_i$ and $\theta_i$. After this training process, DNN can obtain the parameters of the estimated normal distribution for $y$ from decisions $x$ and external information $\theta$. Therefore, embedding its structure within the stochastic optimization problem will allow the possibility of knowing the probabilistic behavior of $y$ conditioned on the decisions about $x$ while decision-making takes place.

\subsection{Distributional Neural Network Embedding}
\label{sec:dnn_emb}

Returning to (\ref{eq:output_1ayer}), we can notice that the input taken by $h(\cdot)$ is a simple linear combination of fixed weights and variables, which does not require additional modifications to be embedded within the optimization problem.

However, we know that $h(\cdot)$ is typically a non-linear activation function. Several works like \cite{anderson2020strong} have studied how to implement trained neural networks as mixed-integer formulations. For example, the employment of the Rectified Linear Unit (ReLU) as $h(\cdot)$ will lead to a max operator that can easily be converted to a set of piece-wise linear constraint employing some additional binary variables, as shown in (\ref{eq:relu}). This ReLU function outputs $v_i^{l}$, defining the maximum value between 0 and the input.

\begin{subequations}\label{eq:relu}
\begin{align}
&v_i^{l} \geq \tilde{v}_i^{l} \\
&v_i^{l} \leq \tilde{v}_i^{l} - M^{low}(1-j_i) \\ 
&v_i^{l} \leq M^{up} j_i
\end{align}
\end{subequations}

\noindent where $\tilde{v}_i^{l}$ takes the value $b_i^l + w_{ij}^{l} v_j^{l-1}$, $M^{low}$ and $M^{up}$ are a negative lower bound and a positive upper bound for the possible values of $\tilde{v}_i^{l}$, respectively, and $j_i$ is a binary variable that takes the value 0 when ReLU outputs 0, and 1 when it outputs $v_i^{l} = \tilde{v}_i^{l}$.

Recent work such as \cite{anderson2020strong} has developed methods to tighten the bound between $M^{low}$ and $M^{up}$, and therefore reduce the computational requirements of the brunch and cut algorithm. However, we consider these procedures to be outside the scope of this work.

Finally, the linear combination of weights jointly with the piece-wise linearized ReLU function will allow us to embed the complete structure of the DNN within the optimization problem (\ref{eq:dcl}). The only difference from a classical NN lies in the output later, where we will obtain an estimation of the mean and standard deviation of the response variable $y$.

\subsection{Scenario generation for learned variables}
\label{sec:sce_gen}

To generate scenarios $\omega \in \Omega$, the different values of $y_{\omega}$ must be sampled from an estimated distribution $\mathcal{Y}(\hat{\mu}_y, \hat{\sigma}_y)$. This process is straight-forward under the normality assumption ($\tilde{y} \sim \mathcal{N}(\hat{\mu}_y, \hat{\sigma}_y)$), which allows to turn the optimization problem (\ref{eq:dcl}) into its Sampled Average Approximation (SAA):

\begin{subequations}\label{eq:dcl_2}
\begin{align}
\underset{x, y_{\omega}}{\min} \;& \sum_{\omega} \pi_{\omega} c(x,y_{\omega};\theta,\xi)  \label{eq:dcl_2_OB}\\
\text{s.t.}&\notag \\
  &g(x,y_{\omega})\leq 0  \label{eq:dcl_othercons2}\\
  &\hat{\mu}_y, \hat{\sigma}_y = f^{\mathcal{D}}(x;\theta,\xi) \label{eq:dcl_modelcons2} \\
  &y_{\omega} = z_{\omega} \times \hat{\sigma}_y + \hat{\mu}_y, \quad \forall \omega \in \Omega \label{eq:dcl_yrandom}
\end{align}
\end{subequations}

\noindent where $\pi_{\omega} = 1/|\Omega|$ represents the probability associated with each scenario $\omega$, with $|\Omega|$ as the total number of scenarios we want to generate. We define stochastic values $y_{\omega}$ as $z_{\omega} \times \hat{\sigma}_y + \hat{\mu}_y$ in constraint (\ref{eq:dcl_yrandom}). In this constraint, $z_{\omega}$ are random values sampled from a standard normal distribution outside the optimization problem. With this procedure of ``distandarization'', $y_{\omega} \; \forall \omega$ are random values sampled from a distribution $\mathcal{N}(\hat{\mu}_y, \hat{\sigma}_y)$. We should emphasize that scenarios $y_{\omega}$ are directly conditioned by decisions $x$, as these modify $\hat{\mu}_y$ and $\hat{\sigma}_y$ through (\ref{eq:dcl_modelcons2}), which alters the definition of the scenario (\ref{eq:dcl_yrandom}). Therefore, the distribution of the random values generated for the learned variable will follow the distribution estimated by the DNN. Furthermore, this problem can be easily extended to characterize each scenario with exogenous sources of uncertainty, such as traditional SAA approaches.

Finally, we also consider this approach appropriate to deal with statistical uncertainty for learned constraints in stochastic mixed-integer optimization problems, as it combines both the high predictive performance of the neural network structure and the possibility of generating scenarios within the optimization model. The employment of distributional forecasting instead of single-point estimations opens a wider range of possibilities for the decision-maker. For example, taking into consideration the possible values of the learned constraint, decision-makers can adapt their actions to behave as risk-averse or risk-neutral when dealing with profits or costs conditioned to this learned response.

\section{Case Study: Virtual Power Plant under price-based demand response}
\label{sec:case_studies}

In this section, we test the validity of the novel stochastic DCL approach by studying its performance if used by a Virtual Power Plant (VPP) to set the optimal tariff for its electricity consumers. We will work under the assumption that the consumers' demand is price-responsive and normally distributed. That is, the price the decision-maker offers to consumers will affect the demand, building a suitable constraint learning model for this relationship. The VPP aims to find the optimal set of hourly prices to offer to the consumers for the next 24 hours while dealing with uncertainty in demand response, the day-ahead market price, and the solar power generation. In particular, the VPP can trade energy in a Day-ahead market to supply its consumers, and can operate a storage system to accommodate possible imbalances between its hourly DERs generation and the net demand. As a last resource to cover imbalances, the VPP can buy energy in a balancing market at comparatively much higher prices.

For this case study (and for all practitioners willing to use DCL methodology), an open-access tool is developed: ``DistCL'' \cite{DistCL2022}. This software allows to train a DNN from a given historical dataset, transform it into a set of piece-wise linear constraints, embed it within a stochastic optimization problem, and generate scenario realizations.

\subsection{Notation}

The notation employed to formulate the optimization problem is described in this subsection for quick reference.

\medskip
\noindent Indices and sets:

\begin{itemize}
    \item[--] $T$: Set of hourly periods within a day, indexed by $t$.
    \item[--] $\Omega$: Set of stochastic scenarios, indexed by $\omega$.
\end{itemize}

\medskip
\noindent First-stage variables:
\begin{itemize}
	\item[--] $P_t$: Price offered to the consumers at time $t$.
	\item[--] $Q_{t}^{DA}$: Quantity of energy bought or sold on the day-ahead market at time $t$.
\end{itemize}

\medskip
\noindent Second-stage variables:
\begin{itemize}
    \item[--] $D_{t,\omega}$: Energy demand at time $t$ and scenario $\omega$.
	\item[--] $Q_{t,\omega}^{Bal}$: Quantity of energy bought in the balance market at time $t$ and scenario $\omega$.
	\item[--] $B_{t,\omega}$: Quantity of energy (status) in the battery at time $t$ and scenario $\omega$.
	\item[--] $\Delta_{t,\omega}^{ch}$: Quantity of energy charged into the battery at time $t$ and scenario $\omega$.
	\item[--] $\Delta_{t,\omega}^{ds}$: Quantity of energy discharged from the battery at time $t$ and scenario $\omega$.
\end{itemize}

\medskip
\noindent Parameters:
\begin{itemize}
	\item[--] $\lambda_{t,\omega}^{DA}$: Day-ahead market price at time $t$ and scenario $\omega$.
	\item[--] $\bar{\lambda}_{t}^{DA}$: Hourly expected day-ahead market price at time $t$, calculated as $\frac{1}{|\Omega|} \sum_{\omega} \lambda_{t,\omega}^{DA}$
	\item[--] $S_{t,\omega}$: Solar power produced at time $t$ and scenario $\omega$.
	\item[--] $\lambda^{Bal}$: Energy price in the balance market.
	\item[--] $R^{ds}$ Discharging ramp capacity for the battery.
	\item[--] $R^{ch}$ Charging ramp capacity for the battery.
	\item[--] $\eta^{ds}:$ Efficiency parameter for discharging.
	\item[--] $\eta^{ch}:$ Efficiency parameter for charging.
	\item[--] $B^{max}$: Maximum battery capacity.
	\item[--] $B^{init}$: Initial quantity of energy in the battery.
	\item[--] $\sigma$: Variability to deviate from the expected day-ahead price allowed to the operator.
	\item[--] $\pi_{\omega}:$ Probability assigned to each scenario $\omega$.
	\item[--] $\theta:$ Contextual variables for estimating demand.
\end{itemize}

\subsection{Conceptual Model}

The VPP stochastic optimization problem is defined as follows.

\begin{subequations}\label{eq:stochastic_model_full}
\begin{align}
\underset{\Phi}{\max} \quad & \sum_{\omega \in \Omega} \pi_{\omega} \sum_{t \in T} \left( P_t D_{t,\omega} - \lambda_{t,\omega}^{DA}Q_t^{DA} - \lambda^{Bal}Q_{t,\omega}^{Bal} \right) \label{eq:OB}\\
\text{s.t.}&\notag \\
  &\hat{\mu}_{D_t}, \hat{\sigma}_{D_t} = f^{\mathcal{D}}(P_t, \theta) \quad \forall t \label{eq:sto_consCL1}\\
  &D_{t,\omega} = z_{\omega} \times\hat{\sigma}_{D_t} + \hat{\mu}_{D_t} \quad \forall t,\omega \label{eq:sto_consCL2}  \\
  &D_{t,\omega} \leq Q_t^{DA} + Q_{t,\omega}^{Bal} + S_{t,\omega} + \eta^{ds}\Delta^{ds}_{t,\omega}-\Delta^{ch}_{t\omega} \quad \forall t,\omega \label{eq:sto_cons2}\\
  &\sum_t P_t \leq \sum_t \bar{\lambda}_{t}^{DA} \label{eq:sto_consconf} \\
   &\bar{\lambda}_{t}^{DA}(1-\sigma) \leq P_t \leq \bar{\lambda}_{t}^{DA}(1+\sigma)  \quad \forall t \label{eq:sto_consconf2} \\
  &0 \leq \Delta^{ds}_{t,\omega} \leq R^{ds}  \quad \forall t,\omega \label{eq:sto_cons3a} \\
    &0 \leq \Delta^{ch}_{t,\omega} \leq R^{ch}  \quad \forall t,\omega \label{eq:sto_cons3b} \\
  &B_{t+1,\omega} = B_{t,\omega} - \Delta^{ds}_{t,\omega}+\eta^{ch}\Delta^{ch}_{t,\omega} \quad \forall t<|T|,\omega \label{eq:sto_cons4} \\
  &0 \leq B_{t,\omega} \leq B^{max} \quad \forall t,\omega \label{eq:sto_cons5} \\
  &B_{1,\omega} = B^{init} \leq \frac{1}{|\Omega|} \sum_{\omega} B_{24,\omega}  \quad \forall \omega \label{eq:sto_cons7} \\
  &Q_{t,\omega}^{Bal} \geq 0 \quad \forall t,\omega \label{eq:sto_cons8}
\end{align}
\end{subequations}

\noindent where $\Phi = \{ P_t, Q_t^{DA}, D_{t,\omega}, Q_{t,\omega}^{Bal}, B_{t,\omega}, \Delta_{t,\omega}^{ch}, \Delta_{t,\omega}^{ds} \}$ represents the set of decision variables.

We aim to maximize the objective (\ref{eq:OB}), defined as the expected revenue from consumer demand minus the expected cost of buying (or selling) energy in the day-ahead market and the cost of buying energy from a balance market. Notice that $Q_{t}^{DA}$ can take negative values, representing the possibility of selling exceeding energy into the market.

Constraints (\ref{eq:sto_consCL1})-(\ref{eq:sto_consCL2}) build the distributional constraint learning framework, i.e., we predict the impact of the price tariff ($P_t$) on the consumers' demand ($D_t$) and characterize the uncertainty associated with this response. Under constraint (\ref{eq:sto_consCL1}), the mean $\hat{\mu}_{D_t}$ and the standard deviation $\hat{\sigma}_{D_t}$ of the demand distribution are estimated, which are learned with a Distributional Neural Network $f(\cdot)$ and are dependent on first-stage decisions $P_t$. Stochastic values $D_{t,\omega}$ are generated with a linear constraint (\ref{eq:sto_consCL2}), where $z_{\omega}$ represents random values sampled from a standard normal distribution outside of the optimization problem. Note that the stochastic values $D_{t,\omega}$ keep their dependence on $P_t$ with the fulfillment of both constraints (\ref{eq:sto_consCL1}) and (\ref{eq:sto_consCL2}). 

On the other hand, constraint (\ref{eq:sto_cons2}) ensures that demand is satisfied with the energy sold / bought in the day-ahead market, the energy from the balance market, the stochastic solar power production and the energy charged or discharged from the battery.

The constraint (\ref{eq:sto_consconf}) establishes that the mean price offered to consumers throughout the day is not greater than the daily mean of the expected day-ahead price. This constraint is set to avoid the VPP to fix extremely high prices to consumers who are supposed to benefit from participating in the VPP. Similarly, the constraint (\ref{eq:sto_consconf2}) limits how much the operator can deviate from its offering price from the expected day-ahead electricity price at each hour.

Finally, the constraints (\ref{eq:sto_cons3a})-(\ref{eq:sto_cons5}) limit the charging and discharging capacity within hours and set the actual value of the battery. Constraint (\ref{eq:sto_cons7}) limits the expected value of the battery charge to be at the end of the day at least the same as at the beginning ($B^{init}$).

\subsection{Methodology}

In this problem, we can find three main sources of uncertainty over the predictions that must be addressed: demand, the day-ahead price, and solar power generation. Therefore, specific data and models are needed in each case.

\subsubsection{Data}

Regarding demand, we simulate a demand response from a real-world study using smart-meter data \cite{kiguchi2021predicting}. We will use consumption data from 3337 households in Tokyo, Japan, measured from 1st July 2017 to 31st December 2018 (18 months). To transform these data into a price-responsive one, we will follow an economic model based on the work presented in \cite{sharma2022estimating}.

Therefore, we assume that the demand is elastic with respect to the price and simulate a price-based demand response as presented in (\ref{eq:dr}).

\begin{equation}
\label{eq:dr}
D^R_t = \bar{D}_t \left[ 1 + \eta^R \left( \frac{P_t^r - P_t }{P_t} \right) \right] 
\end{equation}

The new demand $D^R_t$ at time $t$ will be equal to the original demand $\bar{D}_t$ time one plus an elasticity parameter $\eta^R$ times the increment or decrease of the offered price $P_t$ with respect to the reference price for consumers $P_t^r$.

Regarding the day-ahead price of electricity, we also used data from the Tokyo region, in particular the day-ahead market price, which can be found on the Japan Electric Power Exchange (JPEX\footnotemark[1]{})

\footnotetext[1]{JPEX: http://www.jepx.org/english/}

Finally, in relation to solar power production, data on solar power installation can be simulated with the tool provided in the Photovoltaic Geographical Information System by the European Commission\footnotemark[2]{}. There, we can select a maximum peak power generation in a specific geographical point all over the world. Then, we can obtain hourly simulations of the power we would have generated with that specific setting, and data from several related meteorological variables.

\footnotetext[2]{PGIS: https://re.jrc.ec.europa.eu/pvg\_tools/en/\#HR}

\subsubsection{Neural Network Models and Scenario Generation}

To make probabilistic predictions, we will have to train three different machine learning models: one for the demand (the learned constraint to be embedded in (\ref{eq:sto_consCL1})), one for the day-ahead price, and one for solar power generation.

In order to fit a piece-wise linearizable model for the demand, we will employ the DNN presented in Section \ref{sec:distcl}. Regarding the features to train the model, we will employ calendar variables (hour, day of the week, month, week of the year, and a holiday indicator variable), the price offered to the consumer (we train the model assuming we offer the day-ahead price), the temperature (it is already included in the consumption dataset), and different lags for the demand and temperature at 24, 48, ..., 120 hours. Notice that we do not use, for example, the lags of the previous hour, as in the optimization problem we will have to decide at time $t$ from $t+1$ to $t+24$, and one hour lags will not be available for decision-making. The scenario generation with this demand will be the one exposed in Section \ref{sec:sce_gen}, and represented in (\ref{eq:sto_consCL1}) and (\ref{eq:sto_consCL2}) as sets of constraints.

Regarding the day-ahead price and the solar power generation, these are variables that are not dependent on our decisions and therefore can be estimated outside the decision problem with more complex neural network structures. In this sense, we have decided to generate a probabilistic forecast for both variables using Amazon's Autoregressive Recurrent Network (DeepAR) \cite{salinas2020deepar}. This model uses a encoder structure to capture temporal dependencies and generate normally distributed predictions. Both models for day-ahead price and solar power will be trained to forecast for the next 24 hours. To predict the day-ahead price, calendar variables and 24, 48, ..., 120 hours price lags will be employed. In relation to solar power, generation lags at 24, 48 and 72 hours will be employed besides irradiance, temperature, and wind speed variables with their respective lags from 1 hour to 72 hours. These meteorological variables are obtained using the same tool as solar power. Notice that considering 1 hour lags when we aim to estimate the next 24 hours is not completely realistic. Therefore, a random error has been added to these lags representing a forecasting error on the meteorological variables. Finally, the predictions are repeated multiple times for scenario generation.

\subsubsection{Parameter setting}

In this section, we discuss all the parameter setting and configuration needed for solving the optimization problem. First, machine learning models will be trained with data from July 2017 to October 2018; using November 2018 and December 2018 as validation and test sets, respectively. This test period will also be used to solve the optimization problem.

Regarding the price-based demand response, we set the elasticity $\eta^R$ to $0.2$, and the different reference prices $p_t^r$ to the historical median of the day-ahead market price for each hour $t$. 

In relation to the solar power generation, we chose a peak generation of 15000 kW, achieving a mean daily solar production close to half of the mean daily energy demand. The maximum capacity of the battery $B^{max}$ to save this solar energy will be the same as the peak generation: 15000 kWh. To charge and discharge this battery, we set both ramps $R^{ds}$ and $R^{ch}$ at 5000 kW at each time $t$. For both discharging and charging efficiency parameters, $\eta^{ds}$ and $\eta^{ch}$, a value of $0.9$ was selected and an initial power value of the battery $B^{init}$ of 8000 kWh was established.

For the last parameters of the optimization model, the balance price $\lambda^{Bal}$ was set to 100 yen/kWh, the allowed deviation for the price from the expected day-ahead price is set to 25\% ($\sigma = 0.25$), and the number of scenarios generated $|\Omega|$ will be 300. 

\begin{figure}[!ht]
    \centering
    \includegraphics[width=\textwidth]{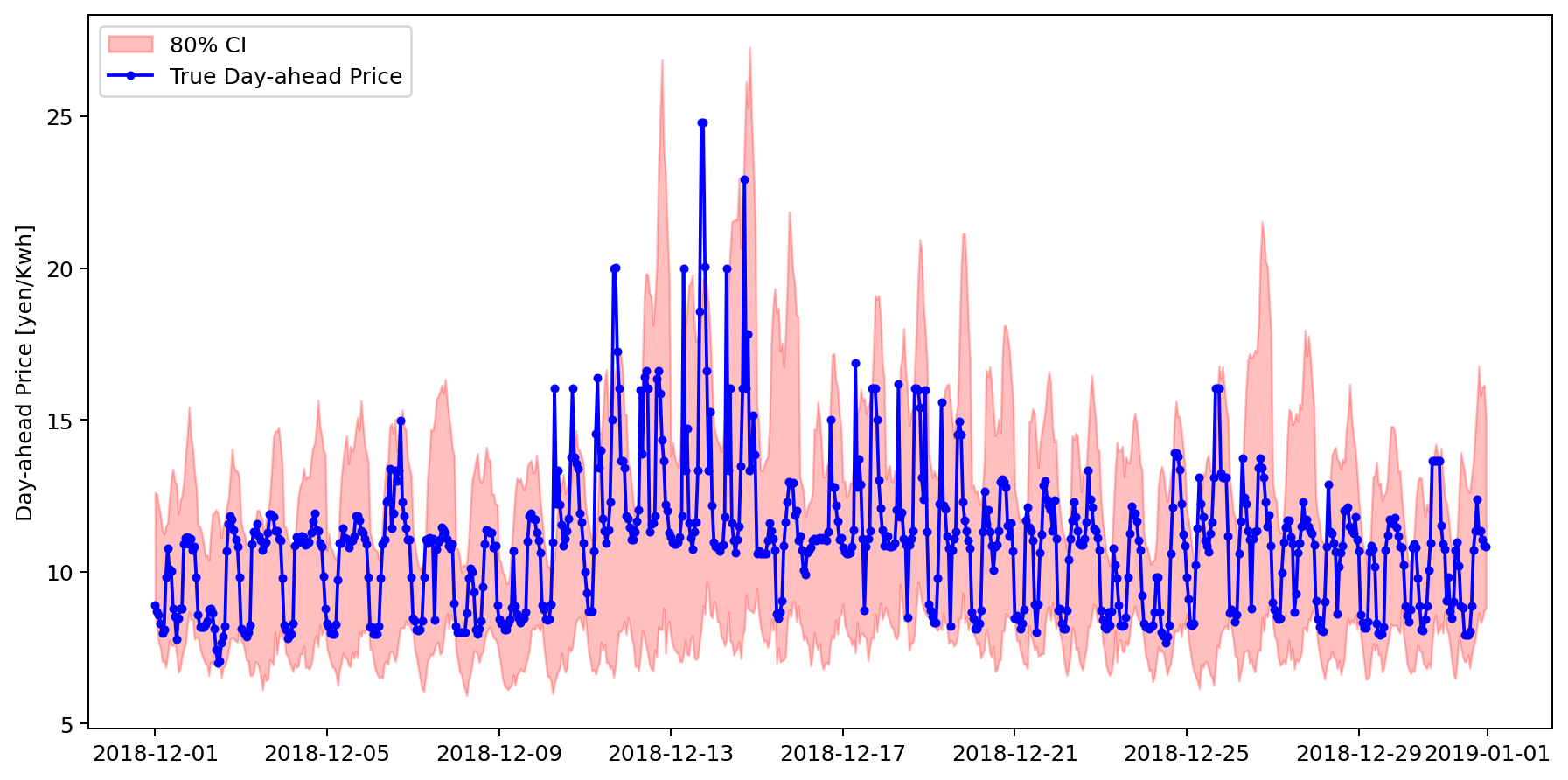}
    \caption{Test probabilistic forecast for day-ahead price}
    \label{fig:da_preds}
\end{figure}

\begin{figure}[!ht]
    \centering
    \includegraphics[width=\textwidth]{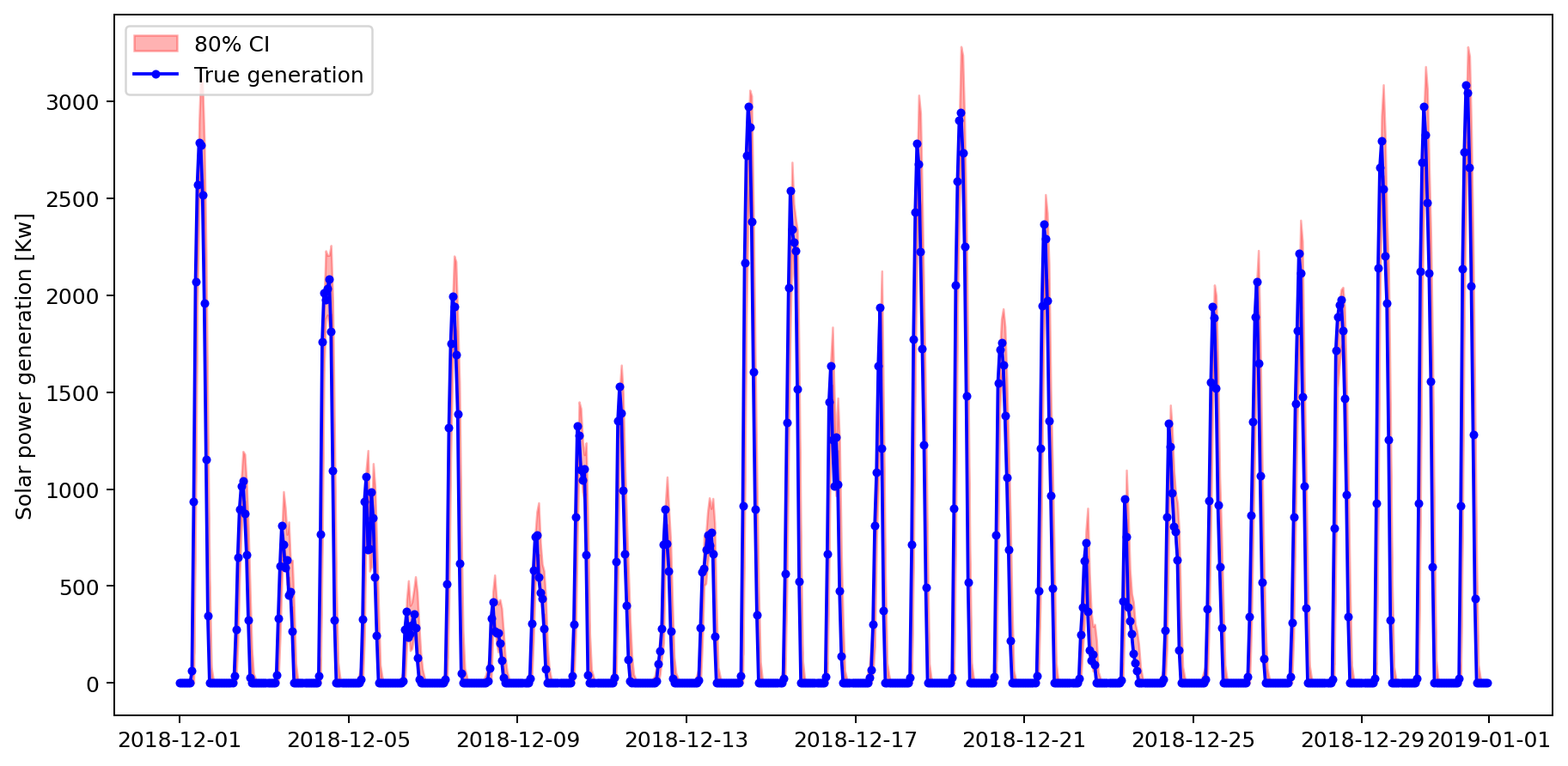}
    \caption{Test probabilistic forecast for solar power generation}
    \label{fig:solar_preds}
\end{figure}

Regarding the DeepAR structure for estimating the day-ahead price and the solar power, a 3 hidden layer with 64 neurons per layer was used for both tasks, employing 120 previous hours for training with respect to the day-ahead price, and 72 hours with respect to the solar power. 300 scenarios are generated from the estimated day-ahead price and solar power production distributions. A visual example of the probabilistic predictions estimated for the price and the solar power can be seen in Figure \ref{fig:da_preds} and Figure \ref{fig:solar_preds}, respectively. Notice that uncertainty is higher when estimating the day-ahead price, as we are trying to estimate a complex response from a market, and we may lack relevant features for training the model.

Finally, in relation to the DNN employed to model the demand, different configurations were used, as will be shown in the following section.

\subsection{Results}

The following results have been obtained using a Python 3.9.12 implementation, using Pyomo 6.3 \cite{hart2017pyomo}. The mathematical solver selected for all computations was Gurobi \cite{gurobi} in its version 9.5. For training the neural network-based models, a Pytorch implementation \cite{NEURIPS2019_9015} was used. Additionally, the computer used included an Intel Core i7 10700 CPU, 64 GB RAM, and a NVIDIA GeForce GTX 2060 graphic card.

\subsubsection{Performance regarding the DNN configuration}

First, we will check the differences in predictive and prescriptive performance when using different DNN configurations to model the demand response with respect to the electricity price-tariff. Let us define $\text{DNN}(l,n)$ as the DNN trained with $l$ hidden layers and $n$ fixed neurons per layer. 

Table \ref{tab:dnn_conf} summarizes the predictive and prescriptive performance of the selected models. The first column indicates the fitting time (seconds needed to reach 5000 epochs). The second shows the Gaussian negative log likelihood loss in the validation dataset. Regarding prescription performance, we indicate first the solving time (in seconds) of the optimization problem (for the complete month of December 2018) and the number of non-binary and binary variables (each day of the month). Then, we show the expected profit during the month of its stochastic distribution and the standard deviation of the sum of profits (one per day in December).

\begin{table}[!ht]
    \centering
    \caption{Predictive and prescriptive performance for different DNN configurations.}
    \label{tab:dnn_conf}
    \resizebox{\textwidth}{!}{%
    \begin{tabular}{c|cc|ccccc}
          & \multicolumn{2}{c|}{Prediction} &  \multicolumn{5}{c}{Prescription} \\ \hline
         Model & Fit. time (s) & Val. loss & Opt. solving time (s) & Non-bin. vars. & Bin. vars. & $\mathbf{E}\left[\sum \text{profit} \right]$ (\yen) & $Var\left[\sum \text{profit} \right]^{1/2}$ (\yen) \\ \hline
         $\text{DNN}(1,20)$ & 598 & -1.597 & 448 & 30445 & 480 & 6034558.66 & 309117.81 \\
         $\text{DNN}(1,50)$ & 633 & -1.730 & 1218 & 31165 & 1200 & 6123409.13 & 279251.66 \\
         $\text{DNN}(2,20)$ & 736 & -1.622 & 644 & 30925 & 960 & 6099529.27 & 284782.92 \\
         $\text{DNN}(2,50)$ & 797 & -1.682 & 1633 & 32365 & 2400 & 6153194.12 & 276791.23 \\
         $\text{DNN}(3,20)$ & 821 & -1.669 & 2352 & 31405 & 1440 & 6108302.05 & 283896.41 \\
         $\text{DNN}(3,50)$ & 926 & -1.706 & 9798 & 33565 & 3600 & 6110752.42 & 281722.50 \\\hline
    \end{tabular}%
    }
\end{table}

The first insight we can obtain by looking at Table \ref{tab:dnn_conf} is that the fitting time increases as the depth of the neural network increases. Besides, the predictive performance of the two- and three-layer DNN is practically the same. The best validation loss is obtained with the model $\text{DNN}(1,50)$.

Regarding the optimization problem solving time, we can notice how it also increases with the size of the DNN, growing exponentially with the size of the last presented model. Note also that the number of binary variables is the same as the number of ReLU operators in the neural network. For example, using model $\text{DNN}(1,20)$ we are dealing with one hidden layer of 20 neurons. In each of these neurons, a ReLU is applied, and therefore 20 binary variables are needed. As there are 24 time periods, we end up with 480 binary variables. 

We can also guess that predictive performance is linked to prescriptive performance. A lower fitting loss leads to a higher expected profit and a lower variance in the results. This is due to the own structure of the predictive model loss, as is composed by the point estimation of the mean (related with an unbiased prediction) and a penalty for high variances. We are interested in fitting an appropriated model with low variance, which entails a less dispersed set of scenarios. Note that in Table \ref{tab:dnn_conf} the model with the lowest validation loss ($\text{DNN}(1,50)$) does not result in the highest expected profit and the lowest variance of the profit ($\text{DNN}(2,50)$). However, this could be due to better performance in the validation set than in the test set (which is unknown at the time of decision-making).

To sum up, we have shown how, in general, predictive performance is related to prescriptive performance, since lower variance is linked to less uncertainty. Furthermore, practitioners should carefully select the size of the DNN, as it could lead to large optimization solving times.

\subsubsection{Comparison with heuristic approaches}
\label{sec:heu}

In this section, we check the real impact of the CL methodology on the result of our optimization problem. In problem (\ref{eq:stochastic_model_full}), allowing the decision-maker to set the price $P_t$ offered to consumers creates the context of CL, as demand is a function of this price and contextual variables. 

For comparison, we consider two heuristic alternatives: The heuristic approach $A$, in which the VPP does not make decisions about the price offered to consumers, and it is fixed to a specific value of 10 yen/kWh; and the heuristic approach $B$, in which the VPP sets the price offered to the consumers externally as the expected day-ahead price, i.e., $\bar{\lambda}_t^{DA}$. From now on, we select $\text{DNN}(1,50)$ as our predictive model. Then, we solve again the optimization problem for the different approaches. Table \ref{tab:heuris} shows the differences regarding the employment of DCL and heuristic approaches.

\begin{table}[!ht]
    \centering
    \caption{DCL and heuristic comparison.}
    \label{tab:heuris}
    \begin{tabular}{c|ccc}
         Approach & Opt. solving time (s) & $\mathbf{E}\left[\sum \text{profit} \right]$ (\yen) & $Var\left[\sum \text{profit} \right]^{1/2}$ (\yen) \\ \hline
            Heuristic A  & 183 & 3104262.74 & 267896.61 \\
          Heuristic B  & 231 & 5383307.74 & 277240.01 \\
          DCL & 1262 & 6123409.13 & 279251.66 \\\hline
    \end{tabular}%
\end{table}

As we can see, when we fix the offered price, the optimization solving time is significantly reduced. Regarding expected profit, when setting a unique offered price (Heuristic $A$) we expect almost half of the profit compared to the DCL approach. Setting the offered price as the expected day-ahead price (Heuristic $A$) rewards a larger expected profit, but around $12\%$ lower than freely deciding on the offered price. In relation to the standard deviation of the sum of profits, it increases as the expected profit increases. However, we do not consider worth-enough to lose $12\%$ of the expected profit to reduce $1\%$ the standard deviation when moving from the complete DCL approach to Heuristic $B$.

\begin{figure}[!ht]
    \centering
    \subfigure[]{\includegraphics[width=0.49\textwidth]{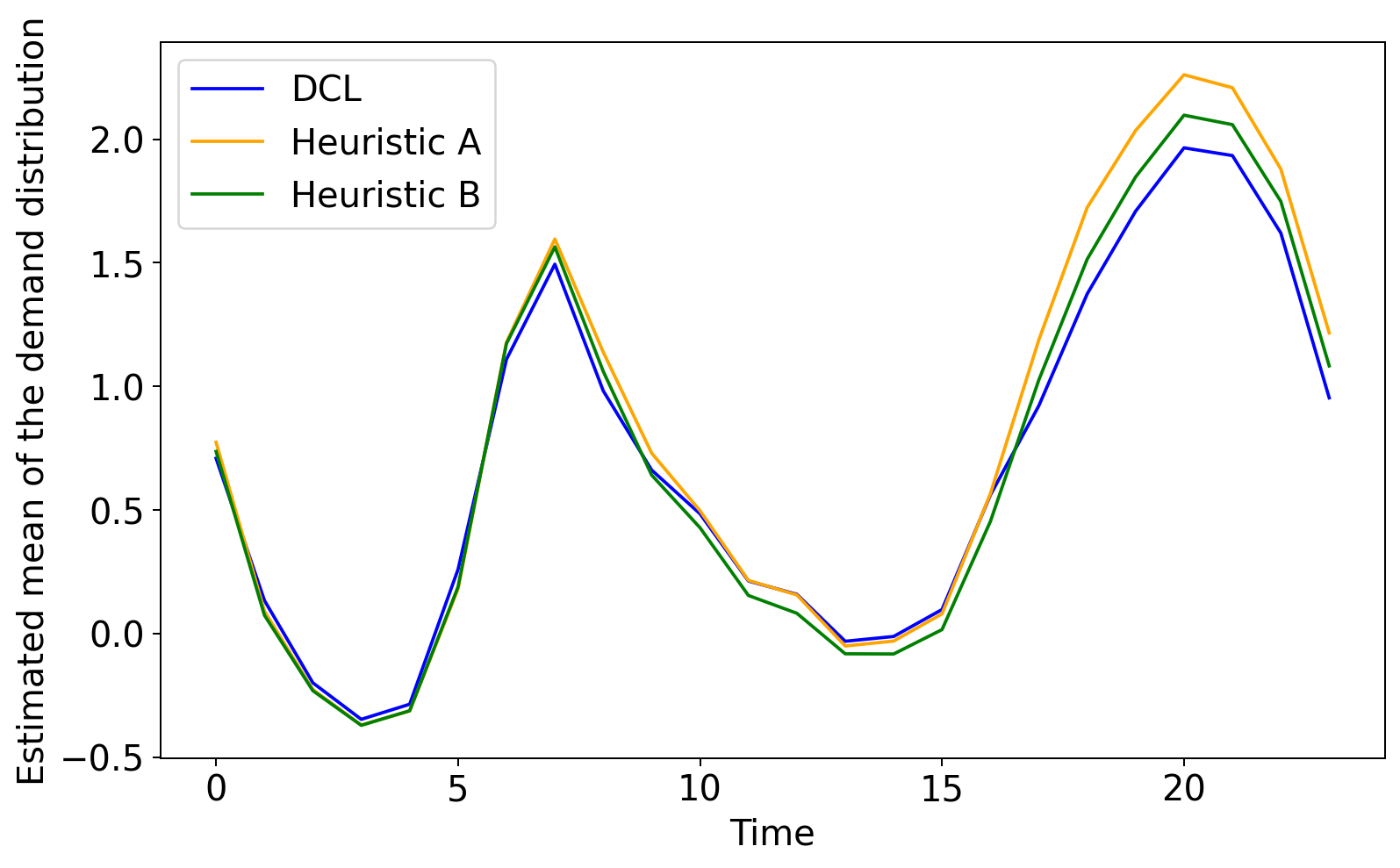}}
    \subfigure[]{\includegraphics[width=0.49\textwidth]{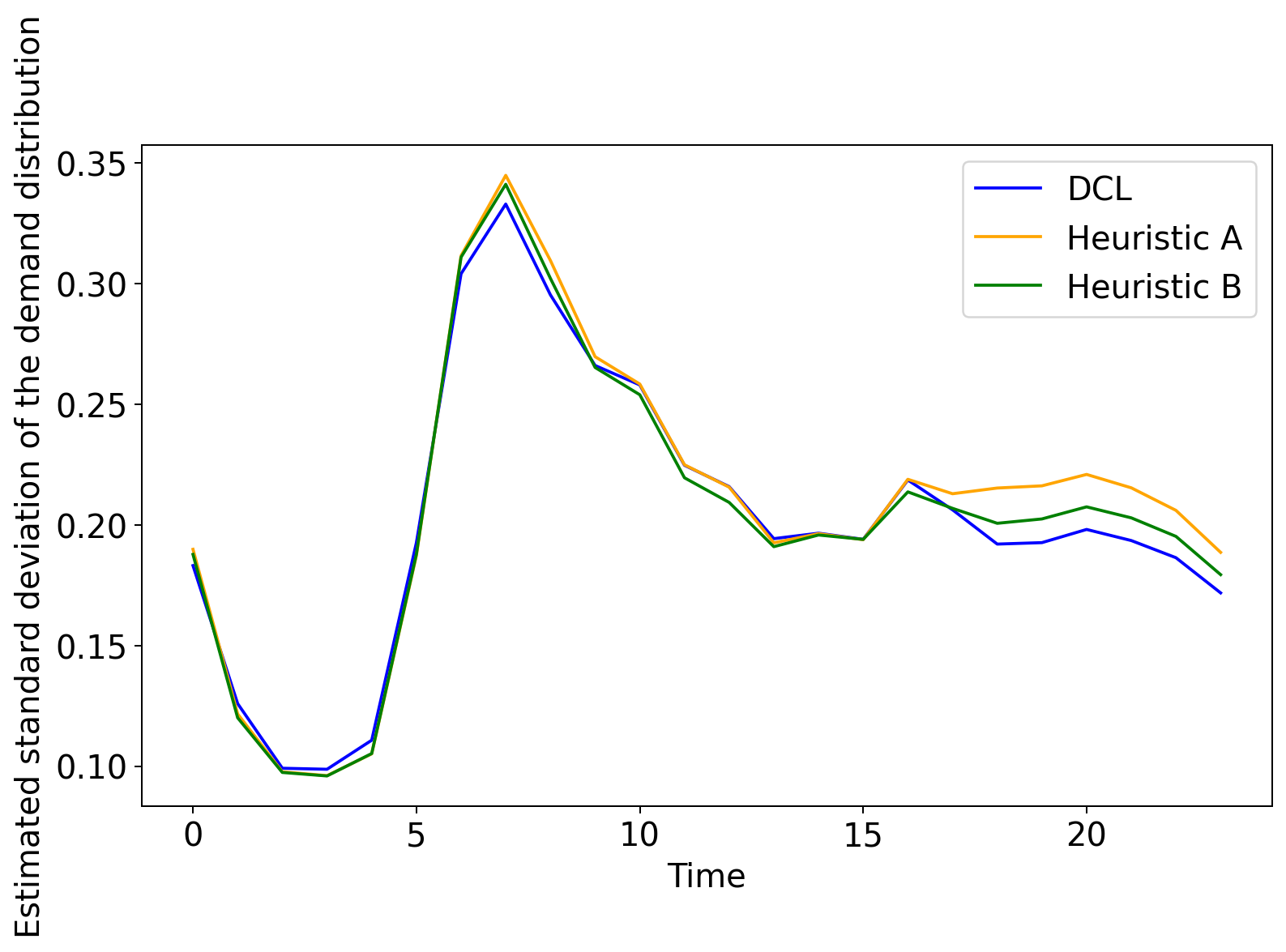}}
    \subfigure[]{\includegraphics[width=0.49\textwidth]{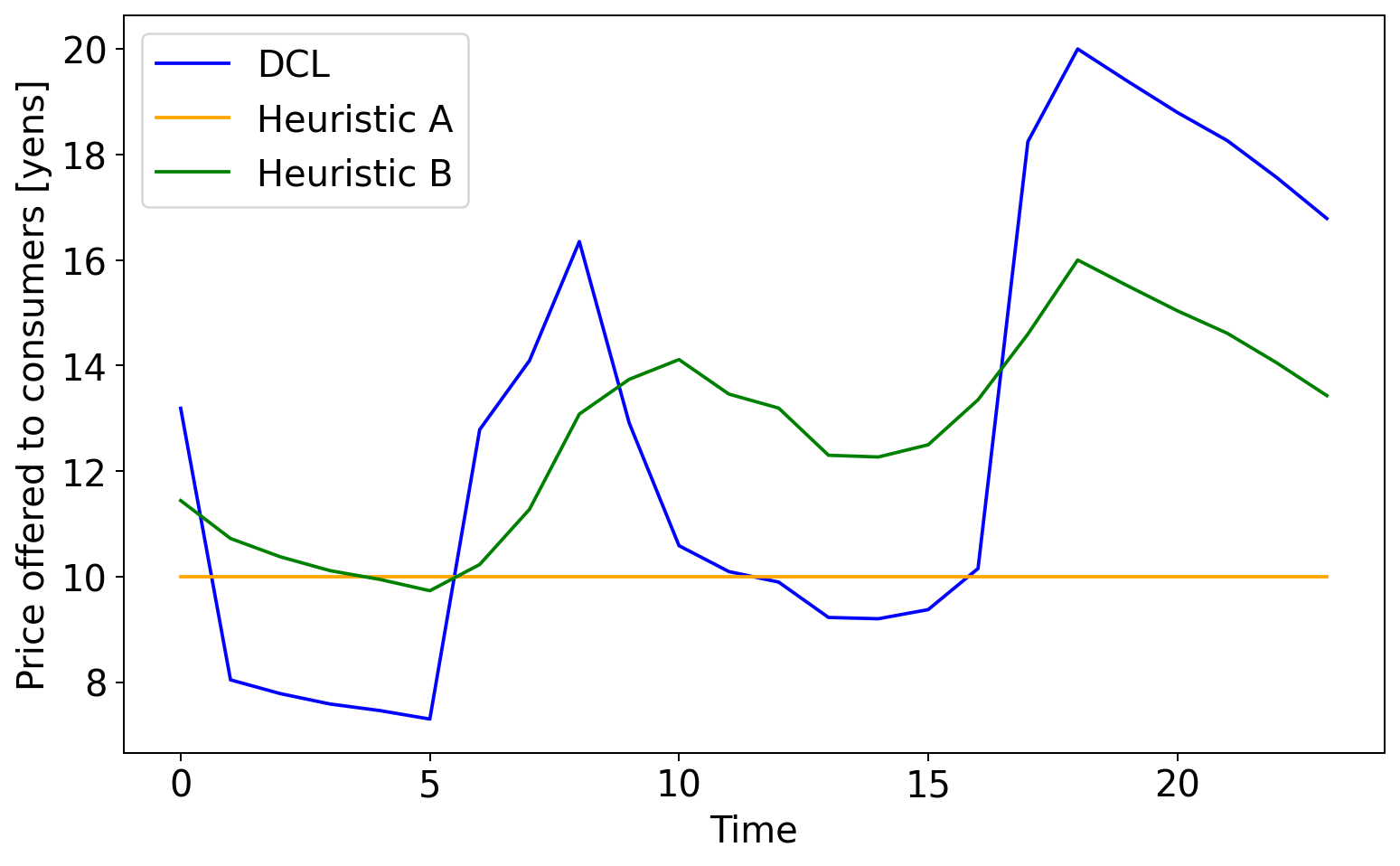}}
    \caption{Estimated mean (a) and standard deviation (b) of the demand distribution. Price offered to consumers (c) (December 12th).}
    \label{fig:demand_heu}
\end{figure}

Finally, we provide a graphical example of the resulting daily demand when employing different approaches, and the decisions made over the price offered to consumers. For December 12th, we obtain the estimated mean and standard deviation of the demand distribution (Figure \ref{fig:demand_heu}(a) and \ref{fig:demand_heu}(b), respectively) when using Heuristic $A$, $B$, or the complete DCL approach. We can also see the price offered to consumers on Figure \ref{fig:demand_heu}(c).

We can notice in Figure \ref{fig:demand_heu}(a) how the decision-maker tries to decrease expected demand during the two peak hours (around 8pm and 8am) without taking into account the valley hours. In terms of the standard deviation of the distribution (Figure \ref{fig:demand_heu}(b)), VPP decisions are able to reduce uncertainty also during the peak hours of demand. Regarding the offering prices (Figure \ref{fig:demand_heu}(c)), the decision-maker further increases the offered price during peak demand hours to reduce consumption peaks.

In this section, we have shown that the CL approach is useful to improve the expected profit of the VPP while maintaining its variance compared with different heuristic approaches. We have also seen how the VPP tries to reduce demand peaks.

\subsubsection{Comparison with deterministic approaches}

In this last section, we want to study how taking into account uncertainty in the constraint learning methodology influences the results in the optimization process compared with point estimation or deterministic approaches, i.e., DCL versus typical CL approaches. To perform this comparison, the following procedure is employed:

\begin{enumerate}
    \item We solve the deterministic version of (\ref{eq:stochastic_model_full}). That is, we solve the problem with a unique scenario, setting the predictive distribution of the day-ahead price and the solar power generation to its mean value. For the CL approach, we set the demand value with the mean estimate of the DNN. This framework will allow us to obtain optimal deterministic prices for consumers ($P_t^{det}$), and the optimal deterministic quantity of energy to sell/buy in the day-ahead market ($Q_t^{DA-det}$). Notice that both variables are considered first-stage decisions.
    \item We solve (\ref{eq:stochastic_model_full}) as a stochastic problem with the complete set of scenarios $\Omega$, but fixing the offering prices $P_t$ to $P_t^{det}$, and the electricity traded on the market $Q_t^{DA}$ to $Q_t^{DA-det}$.
\end{enumerate}

With this approach, we aim to show the differences in the profit distribution when we follow the exposed DCL methodology compared to the distribution obtained when we make first-stage decisions in a non-stochastic context from a deterministic CL methodology.

As an illustrative date, we focus on December 12th, select $\text{DNN}(1,50)$ as our predictive model for the DCL and deterministic methodology, and solve the optimization problem (\ref{eq:stochastic_model_full}).

Table \ref{tab:prof_det} shows the difference in profit distributions for both modeling approaches (DCL vs. Deterministic CL). We show the expected value of the profit distribution, its standard deviation, the Value at Risk (VaR), and the Conditional Value at Risk (CVaR) for $\alpha = 0.1$ \cite{rockafellar2002conditional}.

\begin{table}[!ht]
    \centering
    \caption{Profit distribution for DCL and deterministic approach (December 12th).}
    \label{tab:prof_det}
    \resizebox{\textwidth}{!}{%
    \begin{tabular}{c|cccc}
         Modeling Approach &  $\mathbf{E}\left[ \text{profit} \right]$ (\yen) & $Var\left[\text{profit} \right]^{1/2}$ (\yen) & $VaR_{\alpha=0.1}\left[ \text{profit} \right]$ (\yen)  & $CVaR_{\alpha=0.1}\left[ \text{profit} \right]$ (\yen) \\ \hline
         DCL &  147937.32 & 70301.16 & 58551.76 & 18405.48 \\
          
          Det. CL &  64448.82 & 181838.34 & -132325.70 & -387869.79 \\ \hline
    \end{tabular}%
    }
\end{table}

Comparing DCL with the deterministic approach, we see how the expected value of the profit distribution for the deterministic approach drops dramatically compared to the value from the DCL approach. In the same sense, we can notice that the standard deviation of the distribution for the deterministic approach doubles that of DCL. Another difference from both approaches lies in the VaR and CVaR of the profit distribution, which allows identifying the location of the worst-case profit scenarios. The DCL approach helps to improve CVaR on a larger scale, surpassing the one from the deterministic approach. Solving the optimization problem without considering uncertainty in the first stage implies a very high cost to fulfill the stochastic constraints of demand equilibrium and battery scheduling in the second stage, as the VPP has to deal with some highly unfavorable operation scenarios. As an example, solving the deterministic problem without uncertainty throws an initially estimated profit of almost 160 thousand yens. When this uncertainty is taken into account, the expected profit decreases drastically as shown in Table \ref{tab:prof_det}.

\begin{figure}[!ht]
    \centering
    \includegraphics[width=\textwidth]{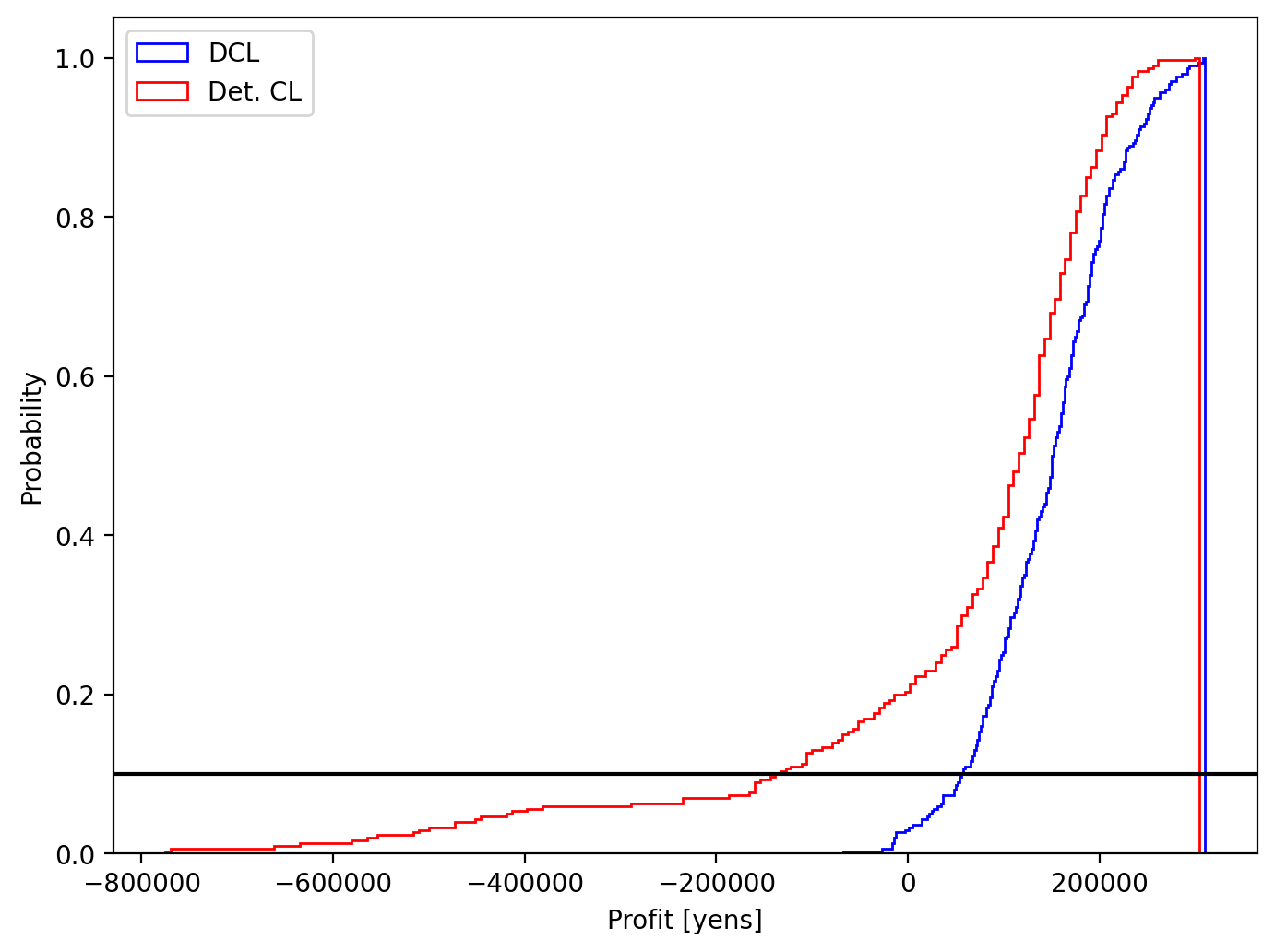}
    \caption{Cumulative Distribution Function of the profit for both DCL and deterministic CL approaches (December 12th).}
    \label{fig:dist}
\end{figure}

We show the differences in the cumulative profit distributions (CDF) in Figure \ref{fig:dist}. As we can see in the profit obtained employing the DCL methodology (red line), only a few values lie below zero, as uncertainty and worst-case scenarios have been taken into account from the beginning. However, when employing a deterministic approach in the first-stage (blue line), around 20\% of the distribution is negative, as these worst scenarios were not taken into account when decision-making takes place. We can also notice how the expected profit, the variance, the VaR and the CVaR are improved when using DCL.

\begin{figure}[!ht]
    \centering
    \subfigure[]{\includegraphics[width=0.49\textwidth]{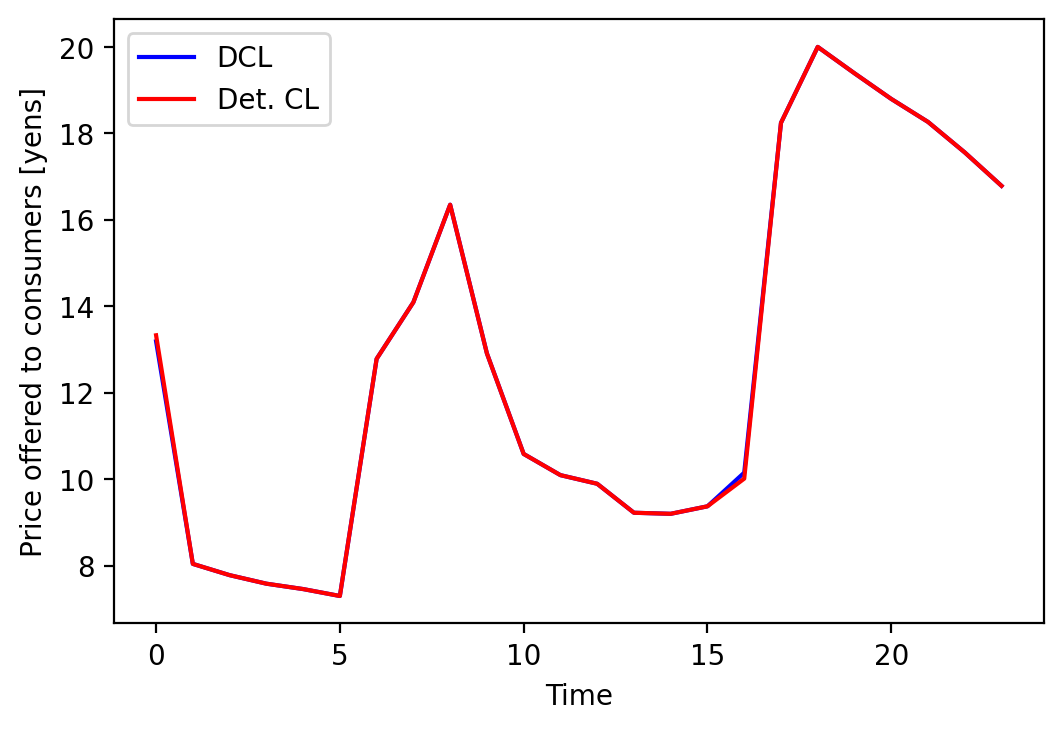}}
    \subfigure[]{\includegraphics[width=0.49\textwidth]{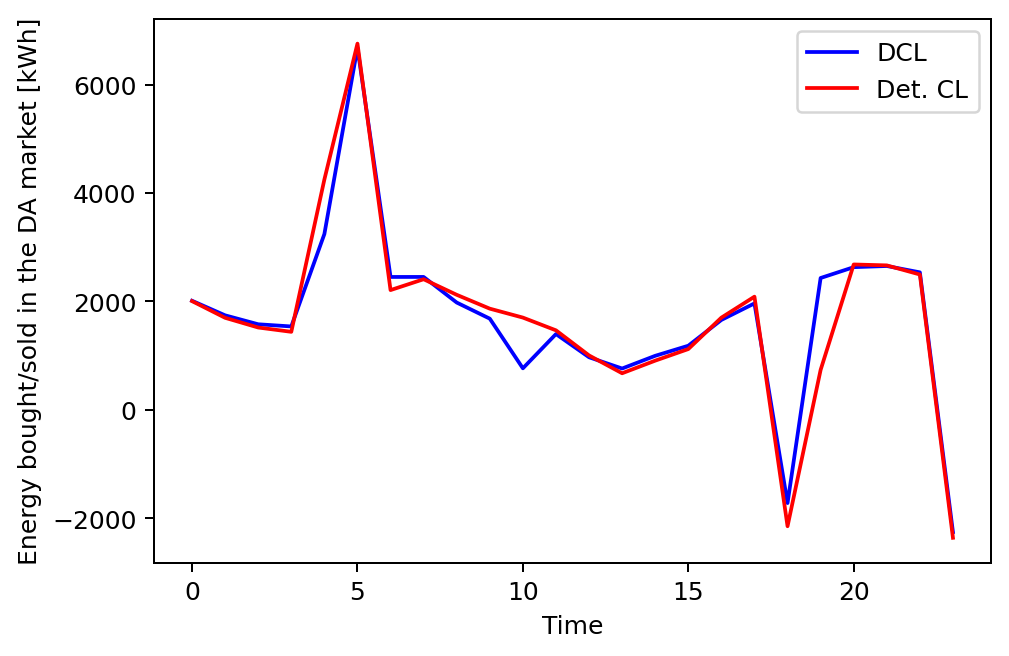}}
    \subfigure[]{\includegraphics[width=0.49\textwidth]{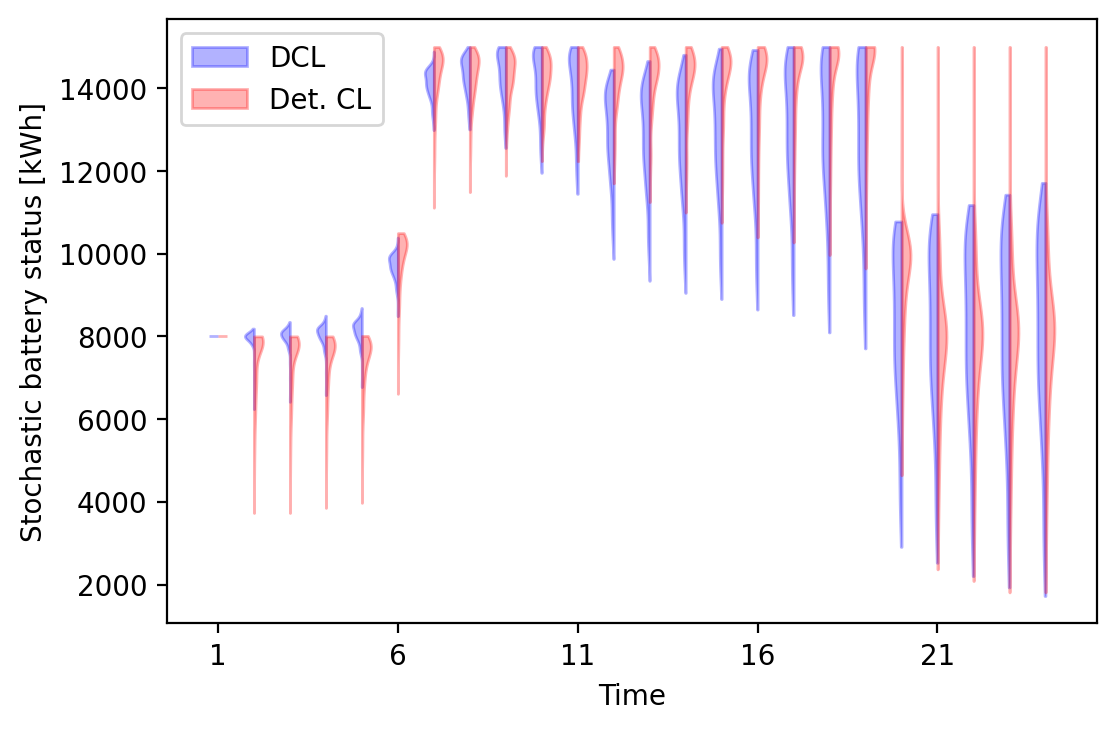}}
    \subfigure[]{\includegraphics[width=0.49\textwidth]{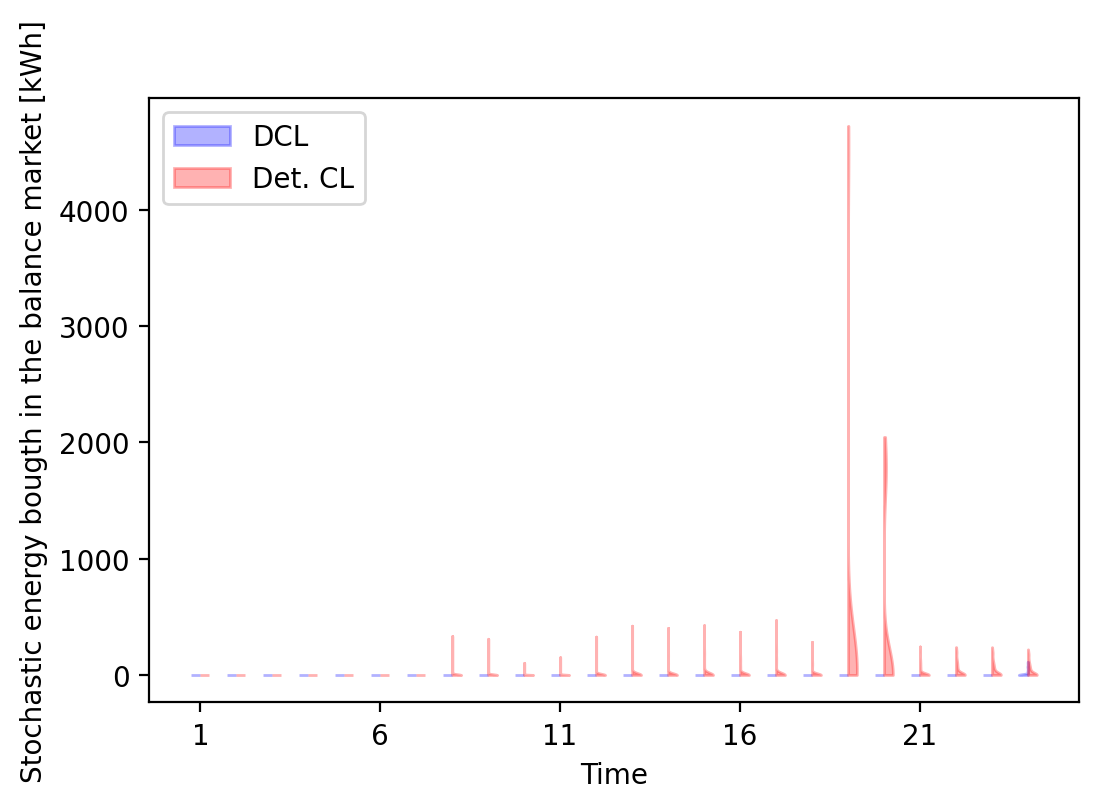}}
    \caption{Offered price to consumers (a), quantity of energy sold/bought in the day-ahead market (b), battery status (c), and energy bought from the balance market (d) employing DCL and deterministic approaches (December 12th).}
    \label{fig:vars_det}
\end{figure}

Finally, we show the value of first-stage and second-stage decision variables for both DCL and the deterministic approach in Figure \ref{fig:vars_det}. Regarding first-stage decision variables, we can check the hourly price offered to consumers in Figure \ref{fig:vars_det}(a), and the energy bought or sold in the day-ahead market in Figure \ref{fig:vars_det}(b). In relation to stochastic second-stage variables, the state of the battery and the amount of energy bought in the balance market are shown in Figure \ref{fig:vars_det}(c) and  \ref{fig:vars_det}(d), respectively.

Regarding the price offered to consumers (Figure \ref{fig:vars_det}(a)), we can notice that it is almost the same for both approaches, as in both cases the VPP tries to reduce demand peaks (see Section \ref{sec:heu} and Figure \ref{fig:demand_heu}). More differences can be noticed in the VPP day-ahead operation (Figure \ref{fig:vars_det}(b)). For instance, the DCL strategy reduces the level of energy sold during the evening hours, probably anticipating scenarios with no energy surplus. These differences in first-stage decision variables, even if they seem small, lead to big changes in the distribution of second-stage variables. For example, in battery status (Figure \ref{fig:vars_det}(c)), the VPP following the DCL approach achieves a lower variance in the distribution, especially in the first and last hours of the day. However, the biggest difference between both approaches can be seen in the energy bought in the balance market (Figure \ref{fig:vars_det}(d)), as the VPP following the DCL approach manages to buy almost no energy, while if the VPP follows the deterministic approach, it must buy large amounts of energy to fulfill the stochastic constraints of demand balance (\ref{eq:sto_cons2}), leading to the decrease in the profit distribution seen in Figure \ref{fig:dist}.

In summary, we have shown in this subsection how taking uncertainty into account is crucial for making optimal decisions, as deterministic approaches lead to lower expected profits and a greater number of worst-case scenarios. For these reasons, we considered that the developed DCL methodology can be of great interest for practitioners to address uncertainty in stochastic contexts and improve the prescriptive performance of standard deterministic CL approaches.

\section{Conclusions}
\label{sec:conclu}

In this work, a novel neural network-based framework has been developed to address uncertainty in the constraint learning methodology in the context of stochastic programming.

State-of-the-art approaches in constraint learning focus on embedding previously trained piece-wise linearizable machine learning models within optimization problems to establish the relationship between a straightforward decision variable $x$ and non-accessible (but $x$-dependent) decision variables $y$. 

However, these machine learning models are fitted to produce a point estimate of $y$ and may not take into account the possible statistical uncertainty about the learned variable $y$, which may be crucial for making optimal decisions in the optimization problem.

For these reasons, we develop the Distributional Constraint Learning (DCL) methodology, which makes use of neural networks to estimate the mean and variance of the conditional distribution of $y$ depending on $x$ and additional contextual information. Assuming a normal distribution over $y$, we are able to embed the neural network structure within the optimization problem and generate scenarios for the estimated distribution, building a stochastic optimization context.

We test the validity of the DCL approach with a Virtual Power Plant (VPP) operation problem, where electricity demands ($y$) can be considered price-responsive ($x$). We show how predictive performance is related to prescriptive one, how the constraint learning approach improves the expected profit over different heuristics, and how DCL outperforms deterministic CL approaches in decision-making under uncertainty.

Future research work could focus on adapting the proposed methodology to determine appropriate uncertainty sets in a robust optimization context. Another extension would be to generalize the normality assumption in the conditional distribution of $y$, and generate scenarios through a set of estimated conditional quantiles.

\section*{Credit authorship contribution statement}

\textbf{Antonio Alc\'antara:} Conceptualization, Methodology, Validation, Investigation, Software, Writing - Original Draft. \textbf{Carlos Ruiz:} Conceptualization, Methodology, Validation, Investigation, Writing - Review \& Editing, Funding acquisition.

\section*{Declaration of competing interest}

The authors declare that they have no known competing financial interests or personal relationships that could have appeared to influence the work reported in this paper.

\section*{Acknowledgements}
The authors gratefully acknowledge the financial support from MCIN/AEI/10.13039/501100011033, project PID2020-116694GB-I00, and the FPU grant (FPU20/00916).


\bibliographystyle{elsarticle-num} 
\bibliography{biblio}

\begin{thebibliography}{10}
\expandafter\ifx\csname url\endcsname\relax
  \def\url#1{\texttt{#1}}\fi
\expandafter\ifx\csname urlprefix\endcsname\relax\def\urlprefix{URL }\fi
\expandafter\ifx\csname href\endcsname\relax
  \def\href#1#2{#2} \def\path#1{#1}\fi

\bibitem{nowotarski2018recent}
J.~Nowotarski, R.~Weron, Recent advances in electricity price forecasting: A
  review of probabilistic forecasting, Renewable and Sustainable Energy Reviews
  81 (2018) 1548--1568.

\bibitem{wan2016probabilistic}
C.~Wan, J.~Lin, Y.~Song, Z.~Xu, G.~Yang, Probabilistic forecasting of
  photovoltaic generation: An efficient statistical approach, IEEE Transactions
  on Power Systems 32~(3) (2016) 2471--2472.

\bibitem{taylor2021combining}
J.~W. Taylor, K.~S. Taylor, Combining probabilistic forecasts of covid-19
  mortality in the united states, European journal of operational research
  (2021).

\bibitem{hao2007quantile}
L.~Hao, D.~Q. Naiman, D.~Q. Naiman, Quantile regression, no. 149 in
  Quantitative Applications in the Social Sciences, Sage, 2007.

\bibitem{khosravi2010lower}
A.~Khosravi, S.~Nahavandi, D.~Creighton, A.~F. Atiya, Lower upper bound
  estimation method for construction of neural network-based prediction
  intervals, IEEE transactions on neural networks 22~(3) (2010) 337--346.

\bibitem{pearce2018high}
T.~Pearce, A.~Brintrup, M.~Zaki, A.~Neely, High-quality prediction intervals
  for deep learning: A distribution-free, ensembled approach, in: International
  conference on machine learning, PMLR, 2018, pp. 4075--4084.

\bibitem{salinas2020deepar}
D.~Salinas, V.~Flunkert, J.~Gasthaus, T.~Januschowski, Deepar: Probabilistic
  forecasting with autoregressive recurrent networks, International Journal of
  Forecasting 36~(3) (2020) 1181--1191.

\bibitem{murty1994operations}
K.~G. Murty, Operations research: deterministic optimization models,
  Prentice-Hall, Inc., 1994.

\bibitem{birge2011introduction}
J.~R. Birge, F.~Louveaux, Introduction to stochastic programming, Springer
  Science \& Business Media, 2011.

\bibitem{bertsimas2020predictive}
D.~Bertsimas, N.~Kallus, From predictive to prescriptive analytics, Management
  Science 66~(3) (2020) 1025--1044.

\bibitem{elmachtoub2022smart}
A.~N. Elmachtoub, P.~Grigas, Smart “predict, then optimize”, Management
  Science 68~(1) (2022) 9--26.

\bibitem{fajemisin2021optimization}
A.~Fajemisin, D.~Maragno, D.~d. Hertog, Optimization with constraint learning:
  A framework and survey, arXiv preprint arXiv:2110.02121 (2021).

\bibitem{kody2022modeling}
A.~Kody, S.~Chevalier, S.~Chatzivasileiadis, D.~Molzahn, Modeling the ac power
  flow equations with optimally compact neural networks: Application to unit
  commitment, Electric Power Systems Research 213 (2022) 108282.

\bibitem{grimstad2019relu}
B.~Grimstad, H.~Andersson, Relu networks as surrogate models in mixed-integer
  linear programs, Computers \& Chemical Engineering 131 (2019) 106580.

\bibitem{maragno2021mixed}
D.~Maragno, H.~Wiberg, D.~Bertsimas, S.~I. Birbil, D.~d. Hertog, A.~Fajemisin,
  Mixed-integer optimization with constraint learning, arXiv preprint
  arXiv:2111.04469 (2021).

\bibitem{mistry2021mixed}
M.~Mistry, D.~Letsios, G.~Krennrich, R.~M. Lee, R.~Misener, Mixed-integer
  convex nonlinear optimization with gradient-boosted trees embedded, INFORMS
  Journal on Computing 33~(3) (2021) 1103--1119.

\bibitem{bergman2022janos}
D.~Bergman, T.~Huang, P.~Brooks, A.~Lodi, A.~U. Raghunathan, Janos: an
  integrated predictive and prescriptive modeling framework, INFORMS Journal on
  Computing 34~(2) (2022) 807--816.

\bibitem{babaei2019data}
S.~Babaei, C.~Zhao, L.~Fan, A data-driven model of virtual power plants in
  day-ahead unit commitment, IEEE Transactions on Power Systems 34~(6) (2019)
  5125--5135.

\bibitem{zhao2021operating}
C.~Zhao, C.~Wan, Y.~Song, Operating reserve quantification using prediction
  intervals of wind power: An integrated probabilistic forecasting and decision
  methodology, IEEE Transactions on Power Systems 36~(4) (2021) 3701--3714.

\bibitem{alcantara2022data}
A.~Alc{\'a}ntara, C.~Ruiz, On data-driven chance constraint learning for
  mixed-integer optimization problems, arXiv preprint arXiv:2207.03844 (2022).

\bibitem{nix1994estimating}
D.~A. Nix, A.~S. Weigend, Estimating the mean and variance of the target
  probability distribution, in: Proceedings of 1994 ieee international
  conference on neural networks (ICNN'94), Vol.~1, IEEE, 1994, pp. 55--60.

\bibitem{anderson2020strong}
R.~Anderson, J.~Huchette, W.~Ma, C.~Tjandraatmadja, J.~P. Vielma, Strong
  mixed-integer programming formulations for trained neural networks,
  Mathematical Programming 183~(1) (2020) 3--39.

\bibitem{DistCL2022}
A.~Alc{\'a}ntara, Distcl, \url{https://github.com/antonioalcantaramata/DistCL}
  (2022).

\bibitem{kiguchi2021predicting}
Y.~Kiguchi, M.~Weeks, R.~Arakawa, Predicting winners and losers under
  time-of-use tariffs using smart meter data, Energy 236 (2021) 121438.

\bibitem{sharma2022estimating}
B.~Sharma, N.~Gupta, K.~Niazi, A.~Swarnkar, Estimating impact of price-based
  demand response in contemporary distribution systems, International Journal
  of Electrical Power \& Energy Systems 135 (2022) 107549.

\bibitem{hart2017pyomo}
W.~E. Hart, C.~D. Laird, J.-P. Watson, D.~L. Woodruff, G.~A. Hackebeil, B.~L.
  Nicholson, J.~D. Siirola, et~al., Pyomo-optimization modeling in python,
  Vol.~67, Springer, 2017.

\bibitem{gurobi}
{Gurobi Optimization, LLC}, \href{https://www.gurobi.com}{{Gurobi Optimizer
  Reference Manual}} (2022).
\newline\urlprefix\url{https://www.gurobi.com}

\bibitem{NEURIPS2019_9015}
A.~Paszke, S.~Gross, F.~Massa, A.~Lerer, J.~Bradbury, G.~Chanan, T.~Killeen,
  Z.~Lin, N.~Gimelshein, L.~Antiga, A.~Desmaison, A.~Kopf, E.~Yang, Z.~DeVito,
  M.~Raison, A.~Tejani, S.~Chilamkurthy, B.~Steiner, L.~Fang, J.~Bai,
  S.~Chintala, Pytorch: An imperative style, high-performance deep learning
  library, in: Advances in Neural Information Processing Systems 32, Curran
  Associates, Inc., 2019, pp. 8024--8035.

\bibitem{rockafellar2002conditional}
R.~T. Rockafellar, S.~Uryasev, Conditional value-at-risk for general loss
  distributions, Journal of banking \& finance 26~(7) (2002) 1443--1471.

\end{thebibliography}

\end{document}